\documentclass[12pt]{article}
\usepackage{multirow}
\usepackage{tabularx}
\usepackage{times}
\usepackage{fullpage,comment}
\usepackage[top=1in, bottom=1in, left=1in, right=1in]{geometry}
\usepackage{amsfonts}
\usepackage{color}
\usepackage{graphicx}
\usepackage[sort,compress]{cite}
\usepackage[utf8]{inputenc}
\usepackage[T1]{fontenc}
\usepackage[mathlines]{lineno}
\usepackage{amssymb,amsthm,amsmath,bm}
\usepackage{algorithmic}
\usepackage{paralist} 
\usepackage{booktabs} 
\usepackage[algo2e,ruled,noend,linesnumbered]{algorithm2e}

\usepackage{authblk} 
\usepackage[colorlinks=true, pdfstartview=FitV, linkcolor=blue,
            citecolor=blue, urlcolor=blue]{hyperref}

\newtheorem{rem}{Remark}[section]
\def\be{\begin{equation}}
\def\ee{\end{equation}}
\def\bes{\begin{equation*}}
\def\ees{\end{equation*}}
\def\bea{\begin{equation} \begin{aligned}}
\def\eea{\end{aligned} \end{equation}}
\def\beas{\begin{equation*} \begin{aligned}}
\def\eeas{\end{aligned} \end{equation*}}
\def\br{\begin{rem}}
\def\er{\end{rem}}
\def\bi{\begin{itemize}}
\def\ei{\end{itemize}}

\def\d{\, \mathrm{d}}

\title{Minimum Reduced-Order Models via Causal Inference}

\author[1]{Nan Chen}
\author[2]{Honghu Liu}
\affil[1]{\footnotesize Department of Mathematics, University of Wisconsin-Madison, Madison, WI 53705, USA (\href{chennan@math.wisc.edu}{chennan@math.wisc.edu})}
\affil[2]{\footnotesize  Department of Mathematics, Virginia Tech, Blacksburg, VA 24061, USA (\href{hhliu@vt.edu}{hhliu@vt.edu})}

\date{{\normalsize \today}}
\begin{document}
\maketitle
\begin{abstract}

Constructing sparse, effective reduced-order models (ROMs) for high-dimensional dynamical data is an active area of research in applied sciences. 
In this work, we study an efficient approach to identifying such sparse ROMs using an information-theoretic indicator called causation entropy. 
Given a feature library of possible building block terms for the sought ROMs, the causation entropy ranks the importance of each term to the dynamics conveyed by the training data before a parameter estimation procedure is performed. It thus allows for an efficient construction of a hierarchy of ROMs with varying degrees of sparsity to effectively handle different tasks. This article examines the ability of the causation entropy to identify skillful sparse ROMs when a relatively high-dimensional ROM is required to emulate the dynamics conveyed by the training dataset. We demonstrate that a Gaussian approximation of the causation entropy still performs exceptionally well even in presence of highly non-Gaussian statistics. Such approximations provide an efficient way to access the otherwise hard to compute causation entropies when the selected feature library contains a large number of candidate functions. Besides recovering long-term statistics, we also demonstrate good performance of the obtained ROMs in recovering unobserved dynamics via data assimilation with partial observations, a test that has not been done before for causation-based ROMs of partial differential equations. The paradigmatic Kuramoto-Sivashinsky equation placed in a chaotic regime with highly skewed, multimodal statistics is utilized for these purposes. 

\medskip
\noindent \textbf{Keywords:} Causation entropy | Data assimilation | Parameter estimation | Kuramoto-Sivashinsky equation | Chaos

\end{abstract}

\tableofcontents

\section{Introduction} \label{Sec_Intro}

Complex dynamical systems appear in many scientific areas, including climate science, geophysics, engineering, neuroscience, plasma physics, and material science \cite{vallis2017atmospheric, strogatz2018nonlinear, wilcox1988multiscale, sheard2009principles,ghil2012topics}. They play a vital role in describing the underlying physics and facilitating the study of many important issues, such as state estimation and forecasting. However, direct numerical simulation is often quite expensive due to the high dimensionality and the multiscale nature of these systems. The situation becomes even more computationally prohibitive when ensemble methods, such as ensemble data assimilation and statistical forecast, are applied to these systems. Therefore, developing appropriate reduced-order models (ROMs) becomes essential not only to reduce the computational cost but also to discover the dominant dynamics of the underlying system.

There exists a vast literature on ROM techniques \cite{ahmed2021closures}. On the one hand, when the complex nonlinear governing equations of the full system are given, one systematic approach to developing ROMs is to project the full governing model to a few leading energetic modes through a Galerkin method. At the core of such Galerkin projections are empirical basis functions constructed with various techniques such as the proper orthogonal decomposition (POD) \cite{HLB96}, the dynamic mode decomposition (DMD) \cite{rowley2009spectral, schmid2010dynamic}, and the principal interaction patterns (PIPs) \cite{hasselmann1988pips, kwasniok1996reduction}. Once the projection is implemented, supplementing the resulting equations with closure terms is often essential to compensate for the truncation error \cite{carlberg2013gnat, noack2011reduced, taira2020modal, CLM20, xie2018data, snyder2022reduced, chorin2015discrete}. On the other hand, many data-driven methods have been developed to learn the dynamics directly from the observed or simulated data associated with the large-scale features or dominant modes of the full systems \cite{brunton2016discovering, ahmed2021closures, chekroun2017data, lin2021data, mou2021data, peherstorfer2015dynamic, hijazi2020data, smarra2018data, chen2018conditional}. The full system does not necessarily need to be known in such a case. It is worth mentioning that there are also many recent developments in non-parametric ROMs or surrogate models, including those resulting from machine learning methods \cite{wan2017reduced, san2018extreme, chattopadhyay2020superparameterization, chen2021bamcafe, pawar2020data, moosavi2015efficient, chattopadhyay2020deep, chattopadhyay2021towards, chekroun2011predicting, mou2023combining, chen2024physics, srinivasan2024turbulence, chen2022conditional}. Many of these developments focus primarily on providing efficient forecast results, while less effort is put into developing a systematic way of quantifying the importance of each constitutive term in the utilized ROMs.  

When the ROMs are given by parametric forms, one would hope that the model is not only skillful in describing the underlying dynamics but also simple enough to facilitate efficient computations. While dynamical models arising from real-world applications usually have a parsimonious structure \cite{temam1997infinite}, data-driven ROMs derived from these full models typically do not inherit such parsimony. On the contrary, as the underlying full system is typically nonlinear, ROMs obtained from projection methods often contain a large number of nonlinear terms. This is because nonlinear interactions among different spatial modes usually cannot be confined to a small subspace spanned by a few spatial modes unless special cancellation properties exist, taking \eqref{Sparse1} below as an example. Similarly, starting with a comprehensive nonlinear model ansatz, applying a standard regression technique to the observed time series typically leads to a high percentage of the terms with non-zero coefficients. Eliminating the terms with weak contributions to dynamics is crucial in balancing the complexity and accuracy of the resulting ROMs. 

For this purpose, it is important to quantify and rank the contribution of each potential constitutive term in a ROM. By the very fabric of chaotic systems, the dynamical contribution of a given term is not simply characterized by the amplitude of its coefficient or even the energy carried by the term. Indeed, removing terms with even a tiny amount of energy can cause a drastic change in the ROM's dynamics, especially when the underlying full model admits sudden transitions such as bursting behaviors \cite{aubry1993preserving, crommelin2004strategies,armbruster1992phase}. Due to these reasons, an intuitive approach of ranking the model's coefficients and removing the terms with small coefficients can fail. Instead, constrained linear regression subject to $\ell_1$-norm regularization, called the least absolute shrinkage and selection operator (LASSO) regression \cite{santosa1986linear, tibshirani1996regression}, has been widely used to discover sparse models from data \cite{brunton2016discovering, boninsegna2018sparse, cortiella2021sparse, schneider2021learning, schaeffer2018extracting,rish2014sparse}. By adjusting the degree of shrinkage (i.e., the level of the $\ell_1$ regularization), the LASSO regression can achieve a varying degree of sparsity since the $\ell_1$ constraint penalizes the absolute values of the model coefficients, pushing some of them to become precisely zero. However, while the importance of each candidate constitutive term in a ROM can be ranked by repeatedly applying the LASSO regression with different degrees of shrinkage, the total computation in such an approach is costly when the total number of candidate terms is large, a situation typically faced when the dimension of the ROM is not small. It is also known that LASSO does not handle severe multicollinearity well \cite{herawati2018regularized}. 

In this paper, we take instead an information-theoretic approach based on the concept of causal inference \cite{almomani2020entropic, fish2021entropic, kim2017causation, almomani2020erfit, elinger2021information, elinger2021causation,sun2014causation} to efficiently construct skillful sparse ROMs. The use of information-theoretic metrics for identifying dynamical models already has a long history, and some pioneer works can be traced back to the 1970s due to H.~Akaike \cite{akaike1973information, akaike1974new}. Over the years, several entropy-based metrics have been studied to quantify information flow or the potential causal influence of one set of data on another, and their scopes of usage span far beyond model identification \cite{lee1998independent,lozano2022information,kleeman2011information,majda2018model}. For our purpose here, we utilize the concept of {\it causation entropy} recalled in \eqref{Causation_Entropy} below, and explore its ability to rank the relative importance of each potential constitutive term in a ROM. 

Initially introduced in \cite{sun2014causation} in the context of network inference, causation entropy extends the concept of {\it transfer entropy} \cite{schreiber2000measuring} to allow for the conditioning on all other candidate terms' influence when the influence of a particular term is considered. Previous studies illustrated that causation entropy, as well as other related variants, are appropriate surrogates for quantifying the dynamical relevance of each candidate term \cite{almomani2020entropic, elinger2021information, chen2023causality, kim2017causation}. Additionally, as demonstrated in \cite{almomani2020entropic}, under certain circumstances such as robustness toward noise and outliers, a causation-based criterion can lead to better sparse ROMs than those constructed from purely cost-function-based minimization techniques including least squares, LASSO, and SINDy (sparse identification of nonlinear dynamics). 

Essentially, causation entropy quantifies the information brought by a given term to the underlying dynamics beyond the information already contained in all the other terms. A larger causation entropy implies that the associated term has a more significant influence/causal effect on the model concerned. Once a feature library of building block functions is chosen as the candidate constitutive terms, the causation entropy associated with each term in the library is computed based on training data using \eqref{Causation_Entropy}; see Sec.~\ref{Subsec:Causality}. Then, different cutoff thresholds can be used to construct ROMs with different complexity. From a ROM identification perspective, causation entropy offers a natural way of ranking the importance of each potential term in the ROM before a parameter estimate step is invoked to learn their corresponding optimal model coefficients (see Step 3 and Step 4 in Sec.~\ref{Sec_CEM_overview}). This separation of the model structure identification (Step 3 in Sec.~\ref{Sec_CEM_overview}) from the parameter estimation (Step 4 in Sec.~\ref{Sec_CEM_overview}) is a key feature of the causal inference framework that distinguishes it from purely cost-function-based minimization techniques. As we demonstrate in Sec.~\ref{Sec:KSE}, it allows for an efficient construction of a hierarchy of ROMs with varying degrees of sparsity to effectively handle different tasks such as recovering long-term statistics and inferring unobserved dynamics via data assimilation. 

Despite the attractive features of causal inference, previous studies on sparse identification using entropy-based metrics focused mainly on either relatively low-dimensional problems or partial differential equation (PDE) models placed in parameter regimes exhibiting chaotic but nearly quasi-periodic dynamics. This article aims to take a modest next step by examining the ability of causation entropy to identify skillful sparse ROMs when a relatively high-dimensional ROM is required to emulate the PDE's dynamics and statistics. By doing so, we also demonstrate that a Gaussian approximation of the causation entropy can still perform exceptionally well even in presence of highly non-Gaussian statistics, which is encouraging since a direct computation of the causation entropy is infeasible when the feature library contains a large number of candidate functions. We also examine the performance of the obtained ROMs for data assimilation, a test that has not been done before for PDEs using causation-based ROMs. We carry out our experiments for the Kuramoto-Sivashinsky equation (KSE) \cite{kuramoto1976persistent, sivashinsky1977nonlinear}, which is a paradigmatic chaotic system with rich dynamical features. 

As shown in Sec.~\ref{Sec:KSE}, the corresponding causation-based ROMs can indeed successfully reproduce both key dynamical features and crucial statistics of the studied KSE model. We also illustrate that to maintain good modeling skills, the level of sparsity achieved depends both on the type of orthogonal basis used in constructing the ROMs and the goals of the ROMs. In that respect, ROMs built from two types of spatial bases are investigated, including the analytical Fourier and data-driven bases built from POD. Due to the particular form of nonlinearity in the KSE, its Galerkin projections under the Fourier basis already exhibit a very sparse structure. This is not the case for Galerkin projections constructed on a POD basis. It is shown in Sec.~\ref{Sec_ROM_eigenbasis} that the causation-based ROMs can recover almost perfectly the sparse structure in the Fourier-Galerkin projections and also reproduce their dynamical and statistical properties, while the ROMs built from the POD basis require far more terms to achieve a similar level of dynamical and statistical performance as reported in Sec.~\ref{Sec_POD_ROM_results}. We also show within the POD setting that when the goal is switched to performing state estimation using data assimilation, we can use a much sparser ROM even if we only observe the time evolution of a few relatively large-scale POD modes; see Sec.~\ref{Sec_DA_results}. 

The rest of the article is organized as follows. We first outline in Sec.~\ref{Sec_CEM} the general procedure to determine the causation-based ROMs, with in particular the core concept of causation entropy and a computationally efficient approximation of this concept recalled in Sec.~\ref{Subsec:Causality}. The usefulness of this ROM framework is then illustrated on the Kuramoto-Sivashinsky in Sec.~\ref{Sec:KSE} in the context of deriving sparse data-driven ROMs to capture either the long-term statistics of the solution or to recover unobserved dynamics using data assimilation. Some additional remarks and potential future directions are then provided in Sec.~\ref{Conclusion}.

\section{A Causation-Based ROM Framework}  \label{Sec_CEM}

\subsection{Overview} \label{Sec_CEM_overview}

With the background and motivations of using causal inference to construct sparse ROMs clarified in the Introduction, we now describe the procedure to determine causation-based ROMs from training datasets. We break it into four steps, with details for the causal inference step and the parameter estimation step provided in Sec.~\ref{Subsec:Causality} and Sec.~\ref{Subsec:ParameterEstimation}, respectively. 

As a starting point, we assume that the training data correspond to the time evolution of a dynamical quantity for which we aim to build a ROM. The training data can come from either observation or the simulation of a finite- or infinite-dimensional full model. For simplicity, we assume that the data are collected at equally spaced time instants, with the time step size denoted by $\Delta t$. Once such a dataset is available, we adopt the following steps to construct a causation-based ROM.

\begin{enumerate}
  \item[Step 1.] {\it Determining the state vector of the ROM via data compression}. Usually, the dimension of the state space for the training data is much higher than the dimension of the sought ROM. In such cases, a data compression procedure is performed to learn a set of empirical basis functions from the training data and then project the training data onto the subspace spanned by the identified basis functions. This data compression procedure is a standard integral part of any projection-based ROMs, and different techniques are available for constructing these basis functions, as reviewed in the Introduction. The output of this process is a set of scalar-valued time series, $\{a_i^j: i = 1,\ldots, n, j = 1, \ldots, N_t\}$, where $n$ denotes the number of empirical basis functions utilized for the projection, and $N_t$ denotes the total number of time instants at which the data are recorded. Viewing $a_i^j$ as the value of a state variable $a_i$ at time $j \Delta t$, the state vector for the sought ROM is then taken to be $\bm{a} = (a_1, \ldots, a_n)^\mathtt{T}$. 

\smallskip  
While the choice of the dimension, $n$, of the ROM depends apparently on the dynamical nature of the modeled quantity as well as the goals of the ROM, one can also draw insights from existing works on rigorous ROM error estimates to assist with this task; see Remark~\ref{Rmk:ROM_steps}. We also note that there are situations where the state vector consists simply of the variables for which the data are provided, taking the paleoclimate proxy records from Greenland ice cores \cite{Boers_al17} and certain physiological time series \cite{weigend1994time} as particular examples. Then, data compression is not needed, and Step 1 is skipped. Namely, in these situations, the reduction aspect is reflected only in reducing the possible number of terms in the ROM's right-hand side without lowering the dimension of the state space.

\item[Step 2.] {\it Constructing a feature library $\mathbb{F}$ of (scalar) candidate functions, from which the constitutive terms of the sought ROM for $\bm{a}$ are selected in Step 3.} Generally, one should use prior knowledge about the possible model behind the training data to inform a judicious choice of the candidate functions. But in case no prior knowledge is available, some typical choices are monomials in terms of the components $a_1, \ldots, a_n$ of $\bm{a}$ up to a certain degree, as well as trigonometric functions or any other common elementary functions in these variables. Time-dependent forcing terms can also be included in the library in case a non-autonomous vector field is expected. At the same time, the computational cost to perform Step 3 below increases when the number of functions in the library increases. Note also that to avoid degeneracy when computing the covariance matrices in Step 3, $\mathbb{F}$ should not contain constant functions. Instead, a constant forcing term can be added afterwards in Step 4 at the stage of parameter estimation; see Sec.~\ref{Subsec:ParameterEstimation}.

  \item[Step 3.] {\it Computing the causation entropy for each function in the feature library $\mathbb{F}$ to determine the ROM's model structure}. This step ranks the importance (in the sense of direct causal relation \cite{sun2014causation}) of all the functions in $\mathbb{F}$ for the time evolution of each derivative $\mathrm{d} a_i/\mathrm{d} t$; see Sec.~\ref{Subsec:Causality} for details. Subsequently, we can select different cutoff thresholds for the computed causation entropies to construct a hierarchy of ROMs with different sparsity levels. This separation of the model structure identification step from the parameter estimation step is a salient feature of the causal inference framework that distinguishes it from purely cost-function-based minimization techniques. 
    
\item [Step 4.] {\it Estimating both the model parameters for the model structure identified in Step 3 and the associated noise amplitude matrix for the model residual, using, e.g., the maximum likelihood method.} See Sec.~\ref{Subsec:ParameterEstimation}. 
  
\end{enumerate}

We end this subsection with the following remark that provides a couple of more comments about the causal inference framework. 

\br \label{Rmk:ROM_steps}

When the training data is noisy, one may need to de-noise the data before performing Step 1 above. Various noise reduction techniques are available for this purpose, including for instance low-pass filtering \cite{kaiser1977data}, singular-spectrum analysis \cite{vautard1992singular}, and total variation based regularization \cite{rudin1992nonlinear}.  

In general, it is difficult to determine the appropriate dimension of the sought ROM in advance. But typically, the state vector identified in Step 1 needs to capture a significant amount of energy contained in the training data. In certain situations, one can also benefit from {\it a priori} ROM error estimates to get an idea about how large the dimension should be; see, e.g., \cite{KV01,koc2022verifiability}. 

Like other data-driven methods, this causation-based approach also requires data for the time derivative of the state vector $\bm{a}$. As a result, the time step $\Delta t$ used in recording the training data should not be too large, unless the time derivative data can be collected through other means rather than by applying a finite difference scheme to the data of $\bm{a}$ obtained in Step 1.  

It is also worth noting that when constructing ROMs for highly chaotic systems, one typically needs to include closure terms to take into account the impact of the orthogonal dynamics not resolved by the ROMs. Additionally, the ROM vector field may also need to respect certain symmetry or energy conservation constraints to ensure stability and accuracy. See Sec.~\ref{Conclusion} for further discussion along these lines. 

\er

\subsection{Determining model structure using causation inference}\label{Subsec:Causality}

We now provide details about how to carry out Step 3 in the previous subsection, once the state vector $\bm{a}=(a_1,\ldots, a_n)^{\mathtt{T}}$ and a feature library $\mathbb{F}$ are identified according to Step 1 and Step 2, respectively. For this purpose, we assume that $\mathbb{F}$ contains a total of $M$ candidate functions, $f_1, \ldots, f_M$, for some $M > 0$: 
\begin{equation}\label{Library}
  \mathbb{F} = \{f_1,\ldots, f_{m-1}, f_m, f_{m+1}, \ldots, f_M\}.
\end{equation}

\paragraph{The definition of causation entropy.}  
The importance of a candidate function $f_m$ in $\mathbb{F}$ to the dynamics of $a_i$ is measured here by the concept of causation entropy. To introduce this concept, let us denote the time derivative of $a_i$ by $\dot{a}_i$. For any $f_m$ in $\mathbb{F}$, let also $\mathbb{F} \backslash {f}_{m}$ be the set of all functions in $\mathbb{F}$ excluding $f_m$. The causation entropy, $C_{f_{m} \rightarrow \dot{a}_i \mid\left[\mathbb{F} \backslash {f}_{m}\right]}$, measures new information provided to $\dot{a}_i$ by $f_m$ in additional to the information already provided to $\dot{a}_i$ by all the other terms in the library $\mathbb{F} \backslash {f}_{m}$. Namely, $C_{f_{m} \rightarrow \dot{a}_i \mid\left[\mathbb{F} \backslash {f}_{m}\right]}$ quantifies to what extent the candidate function $f_m$ contributes to the right-hand side of the equation for $a_i$. Its precise definition is given by \cite{sun2014causation,sun2015causal}:
\begin{equation}\label{Causation_Entropy}
  C_{f_{m} \rightarrow \dot{a}_i \mid\left[\mathbb{F} \backslash {f}_{m}\right]} \stackrel{\mathrm{def}}{=} H(\dot{a}_i | \left[\mathbb{F} \backslash {f}_{m}\right]) - H(\dot{a}_i|\left[\mathbb{F} \backslash {f}_{m} \right], f_m) = H(\dot{a}_i | \left[\mathbb{F} \backslash {f}_{m}\right]) - H(\dot{a}_i|\mathbb{F}),
\end{equation}
where the term $H(\cdot|\cdot)$ is the conditional entropy, which is related to the Shannon entropy $H(\cdot)$ and the joint Shannon entropy $H(\cdot,
\cdot)$ as follows. For two multi-dimensional random variables $\bm{X}$ and $\bm{Y}$ (with the corresponding states being $\bm{x}$ and $\bm{y}$), the following identity (known as the chain rule) holds \cite[Theorem 2.2.1]{cover1999elements}: 
\be \label{Eq_entropy_relation}
H(\bm{Y} | \bm{X}) = H(\bm{X},\bm{Y}) -   H(\bm{X}), 
\ee
where the involved quantities are defined by 
\begin{equation}\label{Entropies}
\begin{split}
  H(\bm{X}) &= -\int_{\bm{x}} p(\bm{x})\log(p(\bm{x}))\d \bm{x},\\
H(\bm{Y} | \bm{X}) &= -\int_{\bm{x}}\int_{\bm{y}} p(\bm{x},\bm{y})\log(p(\bm{y}|\bm{x}))\d \bm{y}\d \bm{x},\\
H(\bm{X},\bm{Y}) &= -\int_{\bm{x}}\int_{\bm{y}} p(\bm{x},\bm{y})\log(p(\bm{x},\bm{y}))\d \bm{y}\d \bm{x},
\end{split}
\end{equation}
with $p(\bm{x})$ being the probability density function (PDF) of $\bm{x}$, $p(\bm{y}|\bm{x})$ the conditional PDF of $\bm{y}$ given $\bm{x}$, and $p(\bm{x},\bm{y})$ the joint PDF of $\bm{x}$ and $\bm{y}$. Regarding the logarithm function $\log()$ involved in \eqref{Entropies}, commonly used bases are $2$, $10$, and Euler's number $e$. For our purpose, the choice of the base is not essential as long as the same base is used for the calculation of all causation entropies since one can convert from base $a$ to base $b$ with the inclusion of a common conversion factor $\log_a b$ \cite[Lemma 2.1.2]{cover1999elements}. This factor is precisely the conversion factor between different units used to measure entropies. For instance, the unit associated with base 2 is called {\it bit} and the one associated with base $e$ called {\it nat}, and $1\, \mathrm{nat} =  \log_2(e)\, \mathrm{bits}$. To fix ideas, we use base 2 in Sec.~\ref{Sec:KSE}. 

On the right-hand side of \eqref{Causation_Entropy}, the difference between the two conditional entropies indicates the information in $\dot{a}_i$ contributed by the specific function $f_m$ given the contributions from all the other functions. Thus, it tells if $f_m$ provides additional information to $\dot{a}_i$ conditioned on the other potential terms in the dynamics. 

Note that the causation entropy $C_{f_{m} \rightarrow \dot{a}_i \mid\left[\mathbb{F} \backslash {f}_{m}\right]}$ actually coincides with the conditional mutual information of $\dot{a}_i$ and $f_{m}$ given $\mathbb{F} \backslash {f}_{m}$, usually denoted by $I(\dot{a}_i; f_{m} | \mathbb{F} \backslash {f}_{m})$, which is always non-negative \cite[Sec.~2.5]{cover1999elements}. Interestingly, even though conditional mutual information \cite{wyner1978definition} is introduced much earlier than causation entropy, its usage for model identification does not seem to be explored until very recently \cite{bhola2023estimating,lozano2022information}. Additionally, by using \eqref{Eq_entropy_relation} in \eqref{Causation_Entropy} (see also the second line in \eqref{Entropy_Gaussians}), we get that $C_{f_{m} \rightarrow \dot{a}_i \mid\left[\mathbb{F} \backslash {f}_{m}\right]} = C_{ \dot{a}_i \rightarrow f_{m} \mid\left[\mathbb{F} \backslash {f}_{m}\right]}$. Namely, given $\mathbb{F} \backslash {f}_{m}$, the new information provided to $\dot{a}_i$ by $f_m$ is the same as the new information  provided to $f_m$ by $\dot{a}_i$.

It is also worthwhile to highlight that the causation entropy in \eqref{Causation_Entropy} is fundamentally different from directly computing the correlation between $\dot{a}_i$ and $f_m$, as the causation entropy also considers the influence of the other library functions. If both $\dot{a}_i$ and $f_m$ are caused by another function $f_{m^\prime}$, then $\dot{a}_i$ and $f_m$  can be highly correlated. Yet, in such a case, the causation entropy $C_{f_{m} \rightarrow \dot{a}_i \mid\left[\mathbb{F} \backslash {f}_{m}\right]}$ will be close to zero as $f_m$ is not the causation of $\dot{a}_i$.

\paragraph{An efficient approximation of causation entropy.} We need to compute the causation entropy $C_{f_{m} \rightarrow \dot{a}_i \mid\left[\mathbb{F} \backslash {f}_{m}\right]}$ for each of the $M$ candidate functions in $\mathbb{F}$ and for each component $a_i$ of the $n$-dimensional state vector $\bm{a}$. Thus, there are in total $nM$ causation entropies to be computed, which can be organized into an $n\times M$ matrix, called the causation entropy matrix. Note that the dimension of $\bm{X}$ in \eqref{Entropies} is either $M-1$ (corresponding to $\mathbb{F} \backslash {f}_{m}$) or $M$ (corresponding to $\mathbb{F}$) in the context of calculating the causation entropy given in \eqref{Causation_Entropy}. This implies that a direct calculation of causation entropies involves both the estimation of high-dimensional PDFs and high-dimensional numerical integrations when the number of library functions $M$ is large, which is known to be computationally challenging. As an alternative, we approximate the joint and marginal distributions involved in \eqref{Entropies} using Gaussians. In such a way, the causation entropy can be approximated as follows \cite{ahmed1989entropy}: 
\begin{equation} \label{Entropy_Gaussians}
\begin{split}
C_{\bm{Z} \rightarrow \bm{X} | \bm{Y}} &=H(\bm{X} | \bm{Y})-H(\bm{X} | \bm{Y}, \bm{Z}) \\
& = H(\bm{X},\bm{Y}) - H(\bm{Y}) - H(\bm{X},\bm{Y},\bm{Z}) + H(\bm{Y},\bm{Z})\\
& \approx \frac{1}{2} \Big(\log(\operatorname{det}(\bm{R}_{\bm{X}\bm{Y}}))-  \log(\operatorname{det}(\bm{R}_{\bm{Y}})) -  \log(\operatorname{det}(\bm{R}_{\bm{X}\bm{Y}\bm{Z}})
 + \log(\operatorname{det}(\bm{R}_{\bm{Y}\bm{Z}}))\Big),
\end{split}
\end{equation}
where $\operatorname{det}(\cdot)$ denotes the determinant of a matrix, $\bm{R}_{\bm{X}\bm{Y}\bm{Z}}$ denotes the covariance matrix of the state variables $(\bm{X},\bm{Y},\bm{Z})$, and the other covariance matrices $\bm{R}_{\bm{Y}}$, $\bm{R}_{\bm{X}\bm{Y}}$ and $\bm{R}_{\bm{Y}\bm{Z}}$ are defined in the same way.

Thus, by assuming that all the involved PDFs follow multivariate Gaussian distributions, the computation of causation entropies boils down to estimating covariance matrices and computing the logarithm of the determinants (log-determinants) of these covariance matrices. This is a much more manageable task when the number of library functions $M$ is too large for other entropy-estimation techniques \cite{darbellay1999estimation,kozachenko1987sample,kraskov2004estimating,schreiber2000measuring} to operate effectively while ensuring accuracy and data efficiency. 

Admittedly, when the concerned data exhibit highly non-Gaussian statistics, the use of Gaussian approximations to compute the associated causation entropy can lead to errors. At the same time, as shown in Sec.~\ref{Sec:KSE}, even though the statistics of the full model's dynamics are highly skewed and multimodal (as revealed in the PDF of the kinetic energy shown in Fig.~\ref{Fig_statistics_eigenbasis}), the Gaussian approximation \eqref{Entropy_Gaussians} still performs exceptionally well as verified using a setting in which the true sparsity structure is known; see Fig.~\ref{Fig_Nonzero_terms_eigenbasis}. In that respect, note also that Gaussian approximations have been widely applied to compute various information measurements in the literature and reasonably accurate results have been reported \cite{majda2018model, tippett2004measuring, kleeman2011information, branicki2012quantifying}. Finally, we would like to point out that even the computation of the log-determinants can be expensive when $M$ is too large, and some further discussion about this is provided in Sec.~\ref{Conclusion}.

\paragraph{Determining the model structure of the ROM.} With the $n\times M$ causation entropy matrix in hand, the next step is determining the model structure. This can be done by setting up a threshold value for the causation entropies and retaining only those candidate functions with the causation entropies exceeding the threshold. When there is a visible gap in the causation entropies (see, e.g., Figure~\ref{Fig_CEM_coef_eigenbasis}), it can serve as a strong indication to set the threshold within this gap. Otherwise, the threshold can be chosen to enforce that a given percentage of the terms in the feature library is kept in the ROM, allowing thus for a hierarchy of ROMs with varying degrees of sparsity by changing the cutoff threshold accordingly; see Sec.~\ref{Sec_POD_ROM_results}. It should also be emphasized that determining the importance of the terms using causation entropy fundamentally differs from that by first ranking the absolute values of the ROM's model coefficients learned from regression and then removing the terms with small coefficients. The latter does not explicitly quantify statistical significance of the eliminated terms and can lead to ROMs with much less accurate results as shown in Sec.~\ref{Sec_DA_results}.

\subsection{Parameter estimation}\label{Subsec:ParameterEstimation}
The final step is to estimate the parameters in the resulting model. For this purpose, we denote the total number of terms in the identified model structure from Step 3 by $s$, and we use a column vector $\bm{\Theta} \in \mathbb{R}^s$ to denote the corresponding $s$ model coefficients to be estimated. We introduce next an $n\times s$ matrix function of the state vector $\bm{a}$, denoted by $\bm{M}$, whose entries are built from the $s$ terms identified in Step 3 as follows. We first put all the $s$ identified terms into each row of $\bm{M}$, arranged in the same order, then for the $i$th row ($i = 1,\ldots, n$), we replace those terms that do not appear in the equation of $a_i$ by $0$. As a result,  the model of $\bm{a}$ can be written in the following vector form: 
\begin{equation}\label{Identified_Model}
  \frac{\d \bm{a}}{\d t} = \boldsymbol\Phi(\bm{a})+\boldsymbol{\sigma} \dot{\bm{W}}(t), \quad \text{with} \quad \boldsymbol\Phi(\bm{a}) = \bm{M}(\bm{a}) \bm{\Theta}.
\end{equation}
In the above model, $\boldsymbol{\sigma} \dot{\bm{W}}(t)$ is a stochastic term, with $\dot{\bm{W}}(t) \in \mathbb{R}^{d\times 1}$ being a white noise for some $d > 0$, and $\boldsymbol{\sigma} \in \mathbb{R}^{n\times d}$ being the noise amplitude matrix. This term $\boldsymbol{\sigma} \dot{\bm{W}}(t)$ aims to model the residual $\frac{\d \bm{a}}{\d t} - \boldsymbol\Phi(\bm{a})$ since usually there does not exist a $\bm{\Theta}$ for which $\boldsymbol\Phi(\bm{a})$ fits perfectly the training data $\frac{\d \bm{a}}{\d t}$ except in some overfitting scenarios. Typically, the dimension of the noise, $d$, is the same as the dimension of the state variable $n$. However, in rare situations when the residual $\frac{\d \bm{a}}{\d t} - \boldsymbol\Phi(\bm{a})$ computed from the training data associated with an estimated $\boldsymbol\Phi(\bm{a})$ has a degenerate covariance matrix, the dimension of $\bm{W}$ would be lower than $n$. 

The parameter vector $\boldsymbol\Theta$ and the noise coefficient matrix $\boldsymbol{\sigma}$ in \eqref{Identified_Model} can be determined using, e.g., the maximum likelihood estimation (MLE) \cite{casella2024statistical}; see \cite{chen2020learning} for the technical details. Notably, the entire parameter estimation can be solved via closed analytic formulae, making the procedure highly efficient. Note also that the optimal parameter $\boldsymbol\Theta$ estimated from the MLE is the same as that obtained from the standard linear regression in the special case when the noise amplitude matrix $\boldsymbol{\sigma}$ is an $n\times n$ diagonal matrix \cite[Chapter 3]{seber2003linear}. MLE can handle constrained minimizations too. Algebraic constraints among certain model parameters can be important depending on specific applications. For instance, the quadratic nonlinearity in many fluid systems represents advection and is a natural energy-conserved quantity. To explicitly enforce energy conservation for the corresponding ROM's quadratic nonlinearity can help prevent the finite time blowup of the solution and also enhance accuracy \cite{majda2012physics, harlim2014ensemble}. Remarkably, closed analytic formulae are still available for parameter estimation in the presence of such constraints; see, e.g., \cite[Section 2.5]{chen2023causality}.

As mentioned in Step 2 of Sec.~\ref{Sec_CEM_overview}, constant functions are excluded in the feature library $\mathbb{F}$ to prevent degeneracy. To add such a constant forcing vector $\bm{b} \in \mathbb{R}^n$ to \eqref{Identified_Model}, we just need to add $n$ additional entries to the parameter vector $\bm{\Theta}$, say after the last entry of the original $\bm{\Theta}$, and also augment $\bm{M}$ with an $n\times n$ identity matrix appended to the right of the last column in the original $\bm{M}$. Usually, if the mean of the residual $\frac{\d \bm{a}}{\d t} - \boldsymbol\Phi(\bm{a})$ computed from the training data is not close to zero, it can be beneficial to add such a constant forcing vector to the ROM.

Note also that if one has precise prior knowledge about a portion of the full model's vector field, it can be advantageous to include it in the ROM. Denoting the projection of this known portion of the vector field onto the ROM subspace by $\bm{Q}$, the drift term $\boldsymbol\Phi(\bm{a})$ in \eqref{Identified_Model} takes then the following form
\begin{equation}\label{RHS_Model}
  \boldsymbol\Phi(\bm{a}) = \bm{M}(\bm{a}) \boldsymbol\Theta + \bm{Q}(\bm{a}),
\end{equation}
where $\bm{Q}\in \mathbb{R}^n$ is a column vector that can depend on $\bm{a}$ but does not involve any free parameters to be estimated. Of course, in this case, we should compute instead the following causation entropies $C_{f_{m} \rightarrow \dot{a}_i \mid (\left[\mathbb{F} \backslash {f}_{m}\right], \bm{Q})}$ for all $i$ and $m$ when determining the model structure in Step 3. 

With the causation-based ROM framework outlined above, we now turn to a concrete application to illustrate its efficiency in identifying effective parsimonious ROMs for the Kuramoto-Sivashinsky equation.

\section{Reduced-Order Models for the Kuramoto-Sivashinsky Equation based on Causal Inference} \label{Sec:KSE}

\subsection{Preliminaries and background}

The Kuramoto-Sivashinsky equation (KSE) \cite{kuramoto1976persistent,sivashinsky1977nonlinear} is a fourth-order dissipative partial differential equation (PDE) that can exhibit intricate spatiotemporal chaotic patterns. It is a prototypical model for long-wave instabilities, which has been derived in various contexts of extended non-equilibrium systems that include unstable drift waves in plasmas \cite{laquey1975nonlinear}, laminar flame fronts \cite{sivashinsky1977nonlinear}, pattern formation in reaction-diffusion systems \cite{kuramoto1976persistent}, and long wave fluctuations in thin film \cite{bertozzi1998long,sivashinsky1980irregular}. Due to its rich dynamical features, the KSE has served as a test ground for various model reduction methods as well as data assimilation techniques in recent years; see, e.g.,  \cite{stinis2004stochastic,lu2017data,CLM20,otto2019linearly,jardak2010comparison,larios2024nonlinear,almomani2020entropic,lunasin2017evolution}.

We consider the one-dimensional KSE:
\be
\label{Eq_KSE}
\partial _t u =  - \nu u_{xxxx} - D u_{xx}   -  \gamma u  u_x,
\ee
which is posed on a bounded interval, $\mathcal{D}=(0,L)$, and subject to periodic boundary conditions. In \eqref{Eq_KSE}, $\nu,D$ and $\gamma$ are positive parameters.

Under the given boundary conditions, since the spatial average is a conserved quantity for the solution $u(x,t)$ of Eq.~\eqref{Eq_KSE}, for simplicity, we restrict to initial data with mean zero by imposing
\be\label{Eq_meanzero}
\int_0^L u(x,t) \d x =0, \quad \forall  \; t \ge 0.
\ee
To build up understanding, throughout the numerical experiments reported below, we will also utilize Galerkin truncations of Eq.~\eqref{Eq_KSE} constructed using either the eigenbasis of the linear operator in \eqref{Eq_KSE} or an empirically constructed basis built from the POD.

\subsubsection{Galerkin projections of the KSE under the Fourier basis}  \label{Sec_Galerkin_eigen}

Due to the assumed periodic boundary conditions, the eigenfunctions of the linear operator $\mathcal{A}=- \nu \partial_{xxxx}  - D \partial_{xx}$ in Eq.~\eqref{Eq_KSE} consist of sine and cosine functions. Thus, the eigenbasis coincides with the Fourier basis. The corresponding Galerkin approximations of Eq.~\eqref{Eq_KSE} can be determined analytically as given below.

First note that the eigenvalues of the linear operator $\mathcal{A}$ subject to the periodic boundary conditions and the additional mean-zero condition \eqref{Eq_meanzero} are given by
\be\label{Eig_KSE}
\beta_n= -\frac{16 \nu \pi^4 n^4}{L^4} + \frac{4D\pi^2 n^2}{L^2}, \quad n \in \mathbb{N},
\ee
where $\mathbb{N}$ denotes the set of all positive integers. Each eigenvalue is associated with two eigenfunctions (labeled by a superscript $\ell$):
\be\label{Modes_KSE}
e_n^{\ell}(x)=\begin{cases}
\sqrt{\frac{2}{L}}\cos\bigg(\frac{2\pi n x}{L}\bigg), \quad \mbox{ if } \ell =0,\\
\sqrt{\frac{2}{L}}\sin\bigg(\frac{2\pi n x}{L}\bigg), \quad \mbox{ if } \ell =1.
\end{cases}
\ee
These eigenfunctions are normalized so that their $L^2(\mathcal{D})$-norm equal to $1$.

Since the eigenfunctions occur in a sine and cosine pair for each wave frequency, we consider Galerkin approximations of Eq.~\eqref{Eq_KSE} of even dimensions. Denote the $2N$-dimensional Galerkin approximation of the KSE solution $u$ under the Fourier basis by
\be \label{Eq_uG_eigen}
u_G(x,t) =  \sum_{n = 1}^N \sum_{\ell = 0}^1 a_n^{\ell}(t) e^{\ell}_n(x).
\ee
Then the amplitudes, $a_n^{\ell}$'s, satisfy the following $2N$-dimensional ODE system
\be\label{KSE_Galerkin_eig}
\frac{\d a_n^\ell} {\d t} = \beta_n a_n^\ell + \sum_{p,q=1}^N \sum_{\ell_1, \ell_2 = 0}^1\Big\langle B(e^{\ell_1}_{p}, e^{\ell_2}_{q}), e_n^{\ell} \Big\rangle a_p^{\ell_1}a_q^{\ell_2}, \quad 1\leq n \leq N, \; \; \ell \in \{0,1\},
\ee
where $B(u,v) =  - \gamma u  v_x$ denotes the quadratic nonlinear term in Eq.~\eqref{Eq_KSE}, and $\langle \cdot, \cdot \rangle$ denotes the $L^2$-inner product for the underlying Hilbert state space.  By direct calculation, we have
\be\label{Sparse1}
\langle B(e^0_p, e^0_q), e^0_n\rangle =
\langle B(e^0_p, e^1_q), e^1_n\rangle =
\langle B(e^1_p, e^0_q), e^1_n\rangle =
\langle B(e^1_p, e^1_q), e^0_n\rangle = 0, \quad \forall \; p,q,n,
\ee
\be\label{Sparse2}
\langle B(e^0_p, e^1_q), e^0_n\rangle
= \langle B(e^1_q, e^0_p), e^0_n\rangle
= \begin{cases}
-\frac{ \gamma \pi n}{\sqrt{2} L^{3/2}}, &\text{ if $n = p+q$},  \vspace{0.2em}\\
\frac{ \gamma \pi (p-q)}{\sqrt{2} L^{3/2}}, &\text{ if $n = |p-q|$},  \vspace{0.2em}\\
0, & \text{otherwise},
\end{cases}
\ee
and
\be\label{Sparse3}
\langle B(e^\ell_p, e^\ell_q), e^1_n\rangle  = \begin{cases}
(-1)^\ell \frac{ \gamma \pi n}{\sqrt{2} L^{3/2}}, &\text{ if $n = p+q$, \; $\ell \in \{0,1\}$}, \vspace{0.2em}\\
\frac{ \gamma \pi n}{\sqrt{2} L^{3/2}}, &\text{ if $n = |p-q|$, \; $\ell \in \{0,1\}$}, \vspace{0.2em} \\
0, & \text{otherwise}.
\end{cases}
\ee
Formulas \eqref{Sparse1}-\eqref{Sparse3} reveal that most of the nonlinear interaction coefficients $\langle B(e^{\ell_1}_{p}, e^{\ell_2}_{q}), e_n^{\ell}\rangle$ in \eqref{KSE_Galerkin_eig} are zero. The resulting Galerkin system \eqref{KSE_Galerkin_eig} has thus a sparse structure. We will show below in Sec.~\ref{Sec_ROM_eigenbasis} that the causal inference criterion presented in Sec.~\ref{Sec_CEM} can be used in a data-driven modeling framework to recover this sparse structure with high fidelity.

\subsubsection{Galerkin projections of the KSE under the POD basis} \label{Sec_Galerkin_POD}

In many applications, empirically computed orthogonal bases can be a more favorable choice than analytic bases due e.g.~to their data-adaptive features. We will thus also assess the skill of the causation inference approach when an empirical basis is used instead. Among the most common choices are the POD \cite{hannachi2007empirical,Holmes_al12,Sir87abc} and its variants \cite{taira2017modal}. Of demonstrated relevance for the reduction of nonlinear PDEs are also the PIPs \cite{hasselmann1988pips,kwasniok1996reduction,kwasniok1997optimal,crommelin2004strategies} that find a compromise between minimizing tendency error with maximizing explained variance in the resolved modes. In the last decade, related promising techniques such as the DMD \cite{rowley2009spectral,schmid2010dynamic,williams2015data} have also emerged; see \cite{tu2013dynamic} for a discussion on the relationships between PIPs, DMD, and the linear inverse modeling \cite{penland1993prediction}.

To fix ideas, we use the POD modes to construct the data-driven Galerkin approximations. Given a cutoff dimension $N$, we denote the POD basis by $\{\varphi_j \; : \; j = 1, \ldots, N\}$, where the basis functions are ranked by their energy content. Recall that the basis functions are orthonormal, i.e., $\langle \varphi_j, \varphi_k \rangle  = \delta_{jk}$ for all $j$ and $k$. The corresponding $N$-dimensional POD-Galerkin system of Eq.~\eqref{Eq_KSE} reads
\be\label{KSE_Galerkin_POD}
\frac{\d a^{\text{POD}}_n} {\d t} = \sum_{j=1}^N A_{nj} a^{\text{POD}}_j + \sum_{i,j=1}^N B_{ij}^n a^{\text{POD}}_i a^{\text{POD}}_j, \quad 1\leq n \leq N,
\ee
where
\be \label{Eq_POD_coef}
A_{nj} =  \langle \mathcal{A} \varphi_j, \varphi_n \rangle, \qquad  B_{ij}^n =  \langle B(\varphi_i,\varphi_j), \varphi_n \rangle,
\ee
with  $\mathcal{A}=- \nu \partial_{xxxx}  - D \partial_{xx}$ and $B(u,v) =  - \gamma u  v_x$ as before.
Once \eqref{KSE_Galerkin_POD} is solved, the corresponding spatiotemporal field that approximates the KSE solution $u$ can be reconstructed via 
\be \label{Eq_uG_POD}
u^{\text{POD}}_G(x,t) =  \sum_{n = 1}^N  a^{\text{POD}}_n(t) \varphi_n(x).
\ee

\subsubsection{Parameter regime and numerical setup} \label{Sec_numerical_setup}

\begin{figure}
\centering
\includegraphics[width=\textwidth]{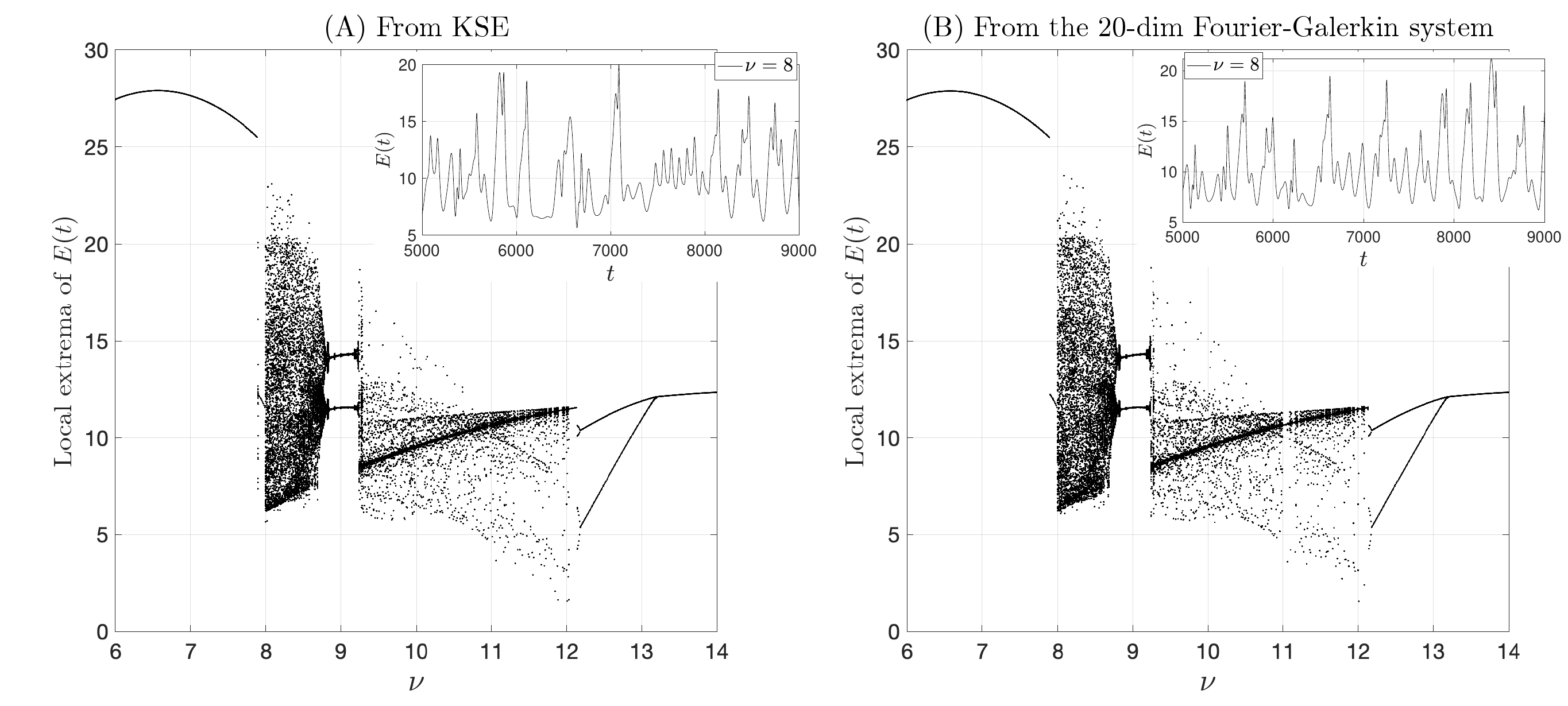}
\caption{\footnotesize Local extrema of the kinetic energy $E$ for the KSE~\eqref{Eq_KSE} (panel A) and its 20-dimensional Fourier-Galerkin approximation (panel B) as the diffusion parameter $\nu$ in \eqref{Eq_KSE} is varied, while the other parameters in \eqref{Eq_KSE} are fixed to be those given in Table~\ref{Table_Regime}. The kinetic energy $E$ for the KSE's solution $u(x,t)$ is defined to be $E(t) = \int_0^L u^2(x,t)\d x$, and for the Galerkin system we have $E(t) = \sum_{n=1}^{10} \sum_{\ell= 0}^1 (a_n^{\ell}(t))^2$, where $a_n^\ell$'s are the state variables of the 20-dimensional Galerkin system of the form \eqref{KSE_Galerkin_eig}. Note that the linear operator $\mathcal{A}=- \nu \partial_{xxxx}  - D \partial_{xx}$ of the KSE has more unstable eigen directions as $\nu$ is decreased; see \eqref{Eig_KSE}. Thus, in principle a regime with a smaller $\nu$ value would require a higher-dimensional Galerkin system to reproduce the KSE's dynamics. We can see that the 20-dimensional Galerkin system recovers vey well the KSE's dynamics for the range of $\nu$ values shown here.}\label{Fig_KSE_bif}
\end{figure}

Throughout Sec.~\ref{Sec:KSE}, we consider the KSE~\eqref{Eq_KSE} in the parameter regime given by Table~\ref{Table_Regime}. For the chosen regime, there are six unstable eigen directions associated with the linear part of the KSE, and the KSE solution is chaotic. As revealed by the ``bifurcation tree'' shown in panel A of Fig.~\ref{Fig_KSE_bif}, the selected regime (with $\nu = 8$) is very close to the borderline in the parameter space where a transition from chaos back to steady state solutions occur as the diffusion coefficient $\nu$ of the stabilizing biharmonic term $-u_{xxxx}$ is further decreased. This interlace between chaotic and non-chaotic dynamics by varying certain model parameter is well documented in the literature, and is not limited to KSE \cite{hyman1986order,CLM17_L9D}. We argue that setting $\nu$ close to this transition borderline creates a challenging scenario to test the performance of the causation-based ROMs, since the discrepancies between the identified ROMs and the full model can easily push the ROM's dynamics into the steady state regime. 

On the implementation side, we solve the KSE by a pseudo-spectral code \cite{kassam2005fourth}, in which the resulting stiff ODE system is integrated using the exponential time-differencing fourth-order Runge-Kutta (ETDRK4) method. We use $N_x=128$ pairs of Fourier modes (hence 128 equally spaced grid points for the spatial domain) to perform the pseudo-spectral discretization, and we use a time step $\Delta t = 0.01$ to integrate the resulting ODE system with the ETDRK4; see Table~\ref{Table_Regime}.

\begin{table}[tbh!]
\caption{System parameters for the KSE \eqref{Eq_KSE}}
\label{Table_Regime}
\centering
\begin{tabular}{cccccc}
\toprule\noalign{\smallskip}
     $\nu$ &  $D$ & $L$  &$\gamma$ & $\Delta t$ & $N_x$    \\
\noalign{\smallskip}\hline\noalign{\smallskip}
 $8$ &  $1$ &  $20\pi$ & 1 & 0.01 & $128$\\
\noalign{\smallskip} \bottomrule
\end{tabular}
\end{table}

It has been checked that a Fourier-Galerkin system with dimension 20 is sufficient to reproduce the dynamical features in the solution of the KSE for the chosen regime; see Fig.~\ref{Fig_KSE_bif}, in which we show that the 20-dimensional Fourier-Galerkin system reproduces accurately the bifurcation tree of the KSE as the diffusion parameter $\nu$ is varied in an interval where the KSE dynamics make transitions among steady states, periodic dynamics, and chaos. While generating a similar bifurcation tree for the 20-dimensional POD-Galerkin system is challenging, since this type of data-driven systems trained for a fixed parameter regime typically cannot be used to infer the dynamics in another parameter regime. Instead, we have checked that the leading Lyapunov exponents for the 20-dimensional Fourier-Galerkin system and the 20-dimensional POD-Galerkin system for the regime given in Table~\ref{Table_Regime} are close to each other, with the leading three exponents for each of these systems given by $0.0087,   0.0066,   0.0049$ and  $0.0088, 0.0067, 0.0048$, respectively. Instead of estimating from generated solution time series, these Lyapunov exponents are computed by evolving the corresponding nonlinear ODE system and their associated variational equations over a long time window [0, 6E5], following \cite[Sec.~3]{wolf1985determining}.\footnote{We did not compute the Lyapunov exponents for the full KSE system due to high computational cost, since the method of \cite[Sec.~3]{wolf1985determining} would require to solve an ODE system of dimension $(128 + 128^2)$ over long time interval.}

The above dynamical insights suggest that the $20$-dimensional Galerkin systems (either constructed from the Fourier basis or the POD basis) approximate reasonably well the full KSE model's dynamics for the chosen parameter regime. We will then test $20$-dimensional causation-based ROMs with different sparsity percentages in the contexts of both Fourier basis and POD basis. We intentionally choose a parameter regime in which the dimension of a high-fidelity Galerkin approximation is not too large in order to not inflate too much the number of candidate functions in the learning library used for computing the causation entropy; see again Sec.~\ref{Subsec:Causality}. See also Sec.~\ref{Conclusion} for some discussions about applying the framework to obtain larger causation-based ROMs with dimensions in the hundreds when needed. 

All the Galerkin approximations of the KSE~\eqref{Eq_KSE}, either constructed from the Fourier basis or the POD basis, are solved using the fourth-order Runge-Kutta method. The causation-based ROMs as well as the thresholded POD-Galerkin systems to be introduced later are simulated with their drift parts approximated using the fourth-order Runge-Kutta method and the stochastic terms approximated using the Euler-Maruyama scheme. The system \eqref{Eq_EnKBF} involved in the data assimilation experiment below is simulated using the Euler-Maruyama scheme for simplicity, which turns out to be sufficient since the involvement of observational data helps alleviate the stiffness of its drift part. The time step size $\Delta t$ for all these models is the same as the one used for solving the PDE itself.

The initial data for the KSE is taken to be $u_0 = \cos(2\pi x/L)$, and the computed solution over the time window $[10^4, 5\times 10^4]$ is used for learning the POD basis as well as for training the related ROMs used in Sec.~\ref{Sec_POD_ROM_results} and Sec.~\ref{Sec_DA_results}. To compute the coefficients involved in the POD-Galerkin system \eqref{KSE_Galerkin_POD}, we approximate each POD basis function using 64 pairs of Fourier modes and then perform the differentiation and integration involved in \eqref{Eq_POD_coef} analytically. For reasons explained in Sec.~\ref{Sec_ROM_eigenbasis}, the causation-based ROMs used in this subsection are trained using the solution for the Fourier-Galerkin systems, still over the time window $[10^4, 5\times 10^4]$.

\subsection{Data-driven inverse models under the Fourier basis} \label{Sec_ROM_eigenbasis}

As pointed out in Sec.~\ref{Sec_Galerkin_eigen}, the Galerkin approximations of the KSE \eqref{Eq_KSE} under the Fourier basis have a sparse structure. Such systems thus provide a good first proof of concept testbed to check whether the causal inference criterion presented in Sec.~\ref{Sec_CEM} can differentiate monomials appearing in the Fourier-Galerkin systems from those that do not when all possible linear and quadratic terms are included in the function library used for model identification.

For this purpose, we place the KSE in the parameter regime given in Sec.~\ref{Sec_numerical_setup} and use the 20-dimensional Fourier-Galerkin system of the form \eqref{KSE_Galerkin_eig} as the true model to generate the training data. We then follow the four-step procedure given in Sec.~\ref{Sec_CEM_overview} to construct the sought ROM, as detailed below. Since we aim to check whether the constructed ROM can recover the sparse structure in the chosen Galerkin system, the state vector of the ROM is the same as that of the Galerkin system. Thus, Step 1 of Sec.~\ref{Sec_CEM_overview} for state vector identification is not needed here. For Step 2, we include all possible linear and quadratic monomials in the feature library $\mathbb{F}$. There are thus a total of 230 candidate functions, consisting of $20$ linear terms and $210$ quadratic terms. To facilitate later discussions, we adopt the following ordering to arrange the 230 library functions. We arrange the 20 unknowns $a_n^\ell$ ($\ell = 0, 1$, $n = 1, \ldots, 10$) of the Galerkin system into a vector $\bm{a} = (a_1^0, \ldots, a_{10}^0, a_1^1, \ldots, a_{10}^1)^{\mathtt{T}}$, and denote the $i$th entry of $\bm{a}$ by $a_i$. The functions in the feature library $\mathbb{F}$ are arranged in the order of
\be \label{Eq_library_order}
\{a_1, \ldots, a_{20}\} \quad \text{ followed by }\quad \{a_ja_k \; |\; 1\le j \le k \le 20\},
\ee
where the following lexicographic order for $(j,k)$ is adopted for the quadratic terms: $(j_1, k_1) < (j_2, k_2)$ if $(j_1 < j_2)$ or $j_1 = j_2$ and $k_1 < k_2$.

\begin{figure}
\centering
\includegraphics[width=\textwidth]{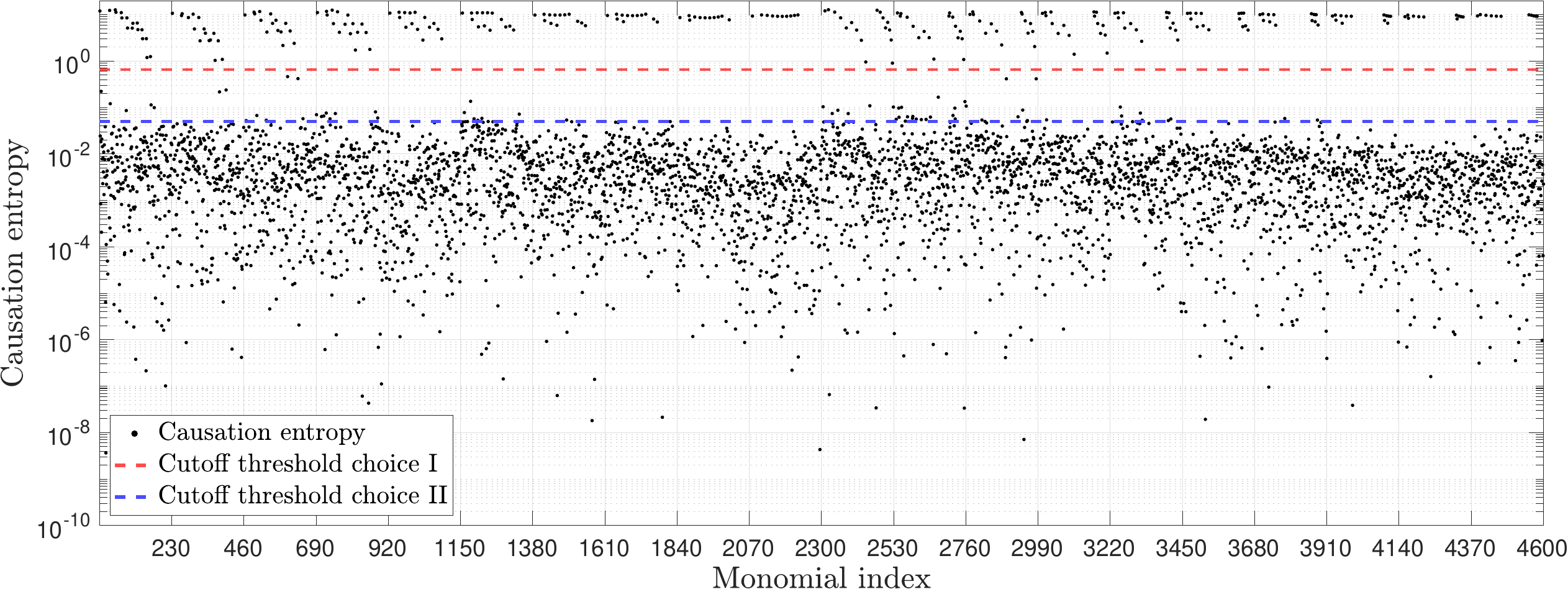}
\caption{\footnotesize Causation entropy that ranks the importance of the library functions for learning a data-driven quadratic inverse model of the 20-dimensional Fourier Galerkin system of the form \eqref{KSE_Galerkin_eig}. The candidate function library $\mathbb{F}$ includes all of the 230 linear and quadratic monomials constructed from the 20 components of the unknown $\bm{Y}$ as listed in \eqref{Eq_library_order}. The causation entropy from each library function $f_m$ to the $i$th equation, $C_{f_{m} \rightarrow \dot{a}_i \mid\left[\mathbb{F} \backslash {f}_{m}\right]}$, is computed according to the approximation formula \eqref{Entropy_Gaussians} given in Sec.~\ref{Subsec:Causality}. The parameter regime is the one given in Sec.~\ref{Sec_numerical_setup}. The causation entropy values are grouped by equation, with the first 230 values (from left) for the first equation, and the next 230 values for the second equation, etc. Also shown are two cutoff thresholds, 0.65 (red line) and 0.05 (blue line). It has been checked that the causation entropy values for all the terms appearing in the true Galerkin system are above the blue line, confirming thus the relevance of this casual inference criterion in identifying constituent terms in the data-driven model.}\label{Fig_CEM_coef_eigenbasis}
\end{figure}

Regarding Step 3 of Sec.~\ref{Sec_CEM_overview} for identifying the ROM's model structure, following the description in Sec.~\ref{Subsec:Causality}, we compute the causation entropy $C_{f_{m} \rightarrow \dot{a}_i \mid\left[\mathbb{F} \backslash {f}_{m}\right]}$ for each library function $f_m$ in $\mathbb{F}$ and each component $a_i$ of $\bm{a}$,  according to the approximation formula \eqref{Entropy_Gaussians}. In total, there are $230\times 20 = 4600$ causation entropy values to compute. These values are shown in Fig.~\ref{Fig_CEM_coef_eigenbasis}, which are grouped by equation, with the first 230 values (from left) for the first equation and the next 230 values for the second equation, etc. There is a visible gap in Fig.~\ref{Fig_CEM_coef_eigenbasis} that separates large causation entropy values (such as those above the red dashed horizontal line) from the smaller ones (below the blue dashed horizontal line), with only very few exceptions falling in between. One is then tempted to suspect that the cutoff threshold for the causation entropy value used for the model structure identification should fall within this gap.

\begin{figure}
\centering
\includegraphics[width=1\textwidth]{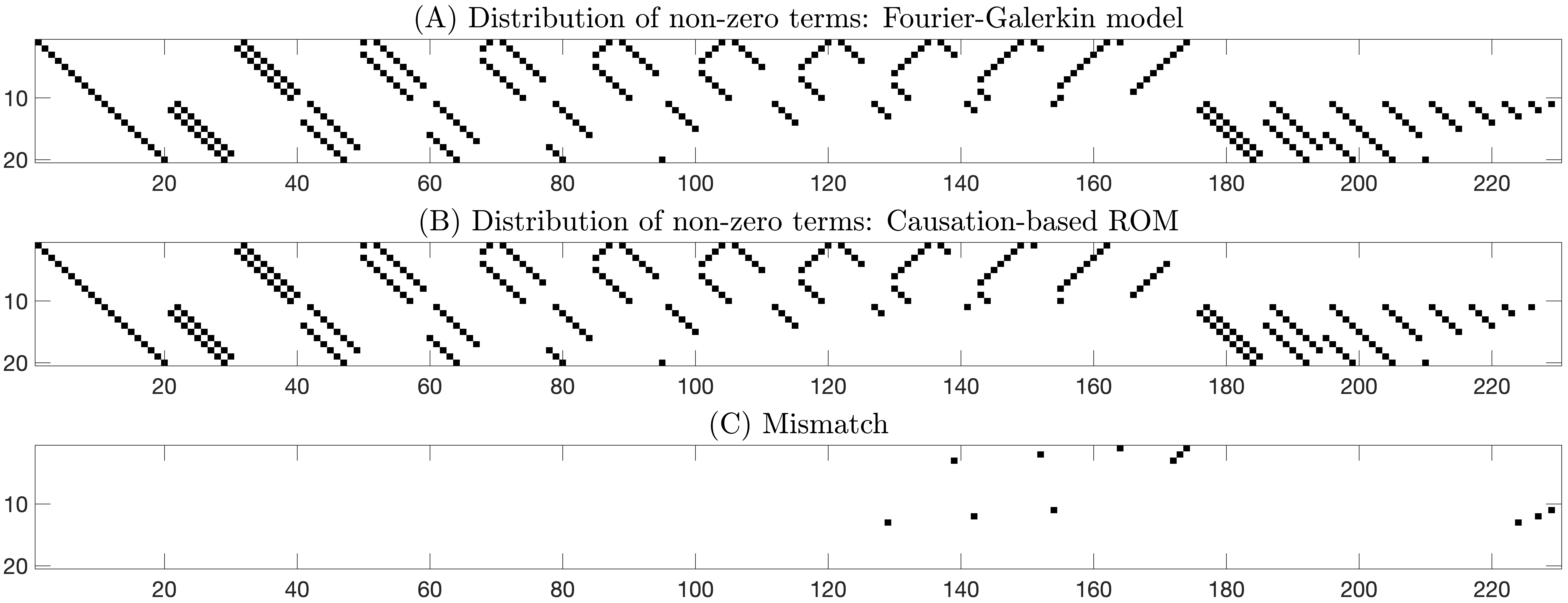}
\caption{\footnotesize Visualization of the distribution of the constituent terms in the 230 function library for the true Galerkin model (panel A) and the learned model (panel B). The learned model is constructed with the cutoff threshold for the causation entropy taken to be 0.65 (red line in Fig.~\ref{Fig_CEM_coef_eigenbasis}). In each panel, the vertical axis consists of 20 rows with each row corresponding to one equation, and the horizontal axis consists of 230 columns with each column corresponding to one function in the learning library. They form thus a $20 \times 230$ mesh. A black square in the $(i,j)$-th grid indicates the $j$th monomial in the library is present in the $i$th equation. The functions in the library are ordered in the way given by \eqref{Eq_library_order}. In particular, the linear terms are placed before the quadratic terms. For instance, the left-most diagonal block in panel A shows that the linear part of the true Galerkin model is diagonal. The mismatches between the learned model and the true model are shown in panel C. There are 12 mismatched terms, all of which are quadratic terms. It has been checked that all these 12 terms are present in the true model, but are missing in the learned model.}
\label{Fig_Nonzero_terms_eigenbasis}
\end{figure}

In the remainder of this section, we use the more severe cutoff threshold marked by the red dashed line in Fig.~\ref{Fig_CEM_coef_eigenbasis} as the cutoff threshold, which corresponds to a numerical value of 0.65 in contrast to 0.05 for the blue dashed line. There are 283 terms whose causation entropy values are above this threshold of 0.65. It turns out all these 283 identified terms are present in the true Galerkin model, which itself has 295 terms. Figure~\ref{Fig_Nonzero_terms_eigenbasis} (panel B) shows the distribution of the identified 283 terms, while that for the true 20-dimensional Galerkin model~\eqref{KSE_Galerkin_eig} is shown in panel A of this figure. The 12 terms in the true model not identified with this cutoff threshold are shown in panel C of Fig.~\ref{Fig_Nonzero_terms_eigenbasis}.

As can be seen in Fig.~\ref{Fig_Nonzero_terms_eigenbasis}, the causation entropy criterion is remarkably successful in identifying the true sparse model. In particular, it has a negligible mismatch rate of $12/(20\times 230) = 0.26\%$ for the total 4600 possible terms to be sifted through. It turns out that all these 12 mismatch terms are quadratic terms, while all the linear terms of the model are correctly identified. Indeed, the linear part of the true model is a diagonal matrix with eigenvalues on the diagonal (see Eq.~\eqref{KSE_Galerkin_eig}). This is represented by the left-most diagonal block in panel A of Fig.~\ref{Fig_Nonzero_terms_eigenbasis} due to the way that the library functions are arranged; see again \eqref{Eq_library_order}. These linear terms are fully captured using the chosen cutoff threshold as shown in panel B of this figure.

\begin{figure}[tbh!]
\centering
\includegraphics[width=1\textwidth]{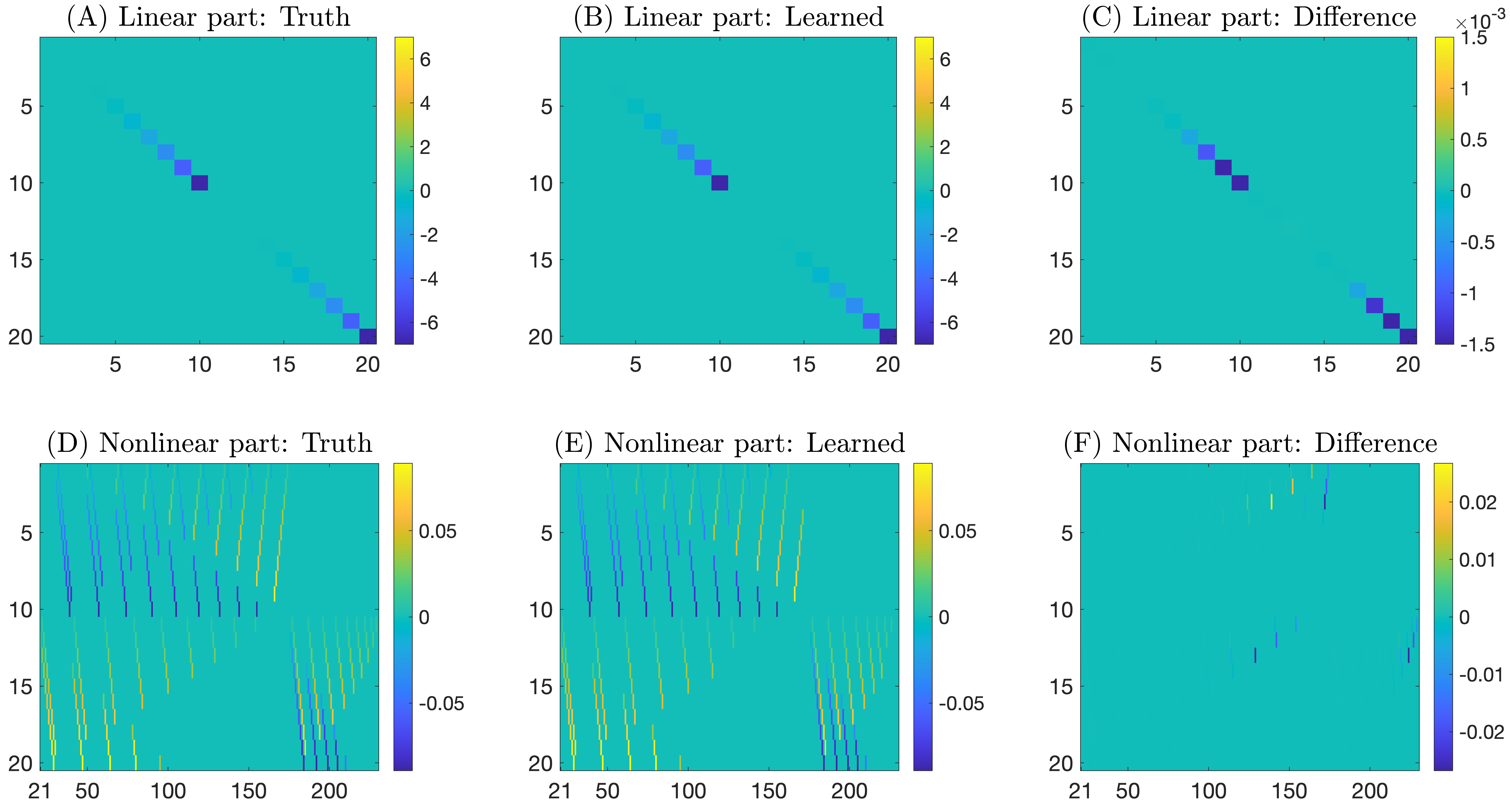}
\caption{\footnotesize Visualization of the model coefficients for the 20-dimensional Fourier-Galerkin system (true model) and the causation-based ROM (learned model). For the true model, all non-zero coefficients occur for the terms marked by black squares shown in panel A of Fig.~\ref{Fig_Nonzero_terms_eigenbasis}. We separated the linear terms (panel A here) from the nonlinear terms (panel D here) for a better visualization, since the coefficients for some of the linear terms are two-order larger than those for the nonlinear terms. The sparsity structure for the learned model is the one shown in panel B of Fig.~\ref{Fig_Nonzero_terms_eigenbasis}.}
\label{Fig_learned_coef_eigenbasis}
\end{figure}

\begin{figure}[tbh!]
\centering
\includegraphics[width=1\textwidth]{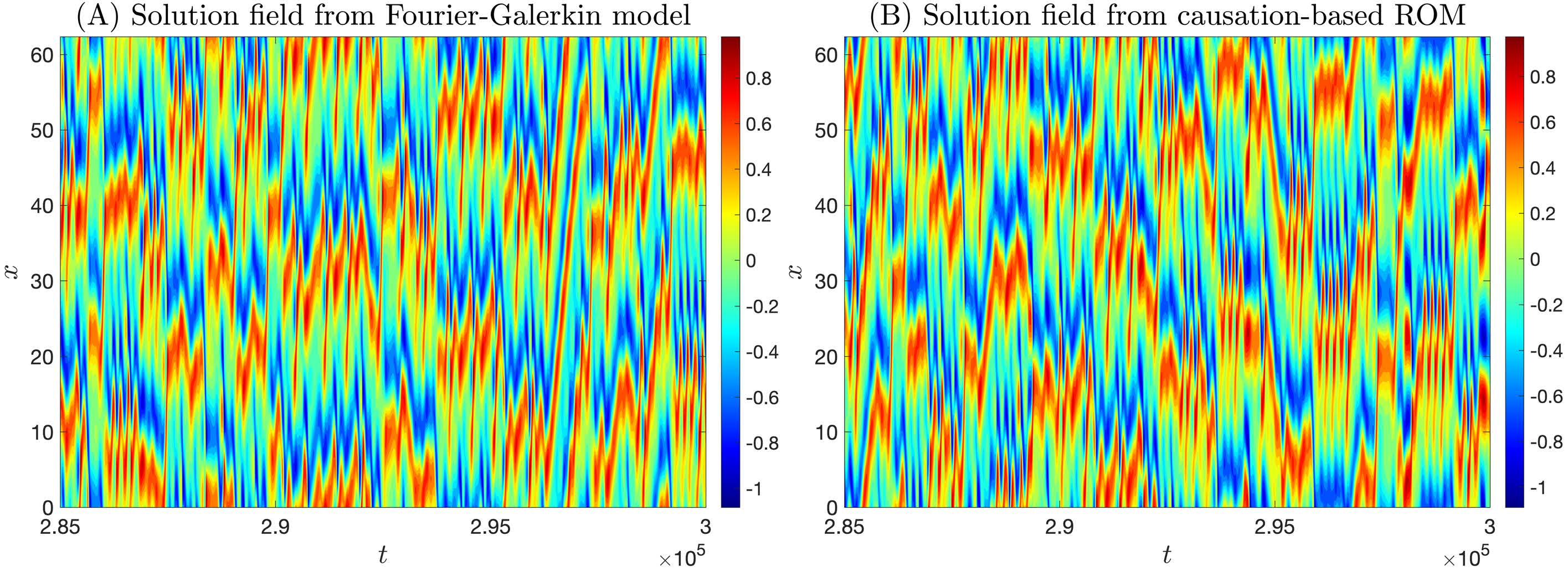}
\caption{\footnotesize Comparison of the reconstructed solution fields obtained using \eqref{Eq_uG_eigen} for the 20-dimensional Fourier-Galerkin system \eqref{KSE_Galerkin_eig} (left panel) and for the causation-based ROM of the same dimension (right panel). The results are shown here in a time window well beyond the training window $[10^4, 5\times 10^4]$.}
\label{Fig_soln_eigenbasis}
\end{figure}

\begin{figure}[tbh!]
\centering
\includegraphics[width=1\textwidth]{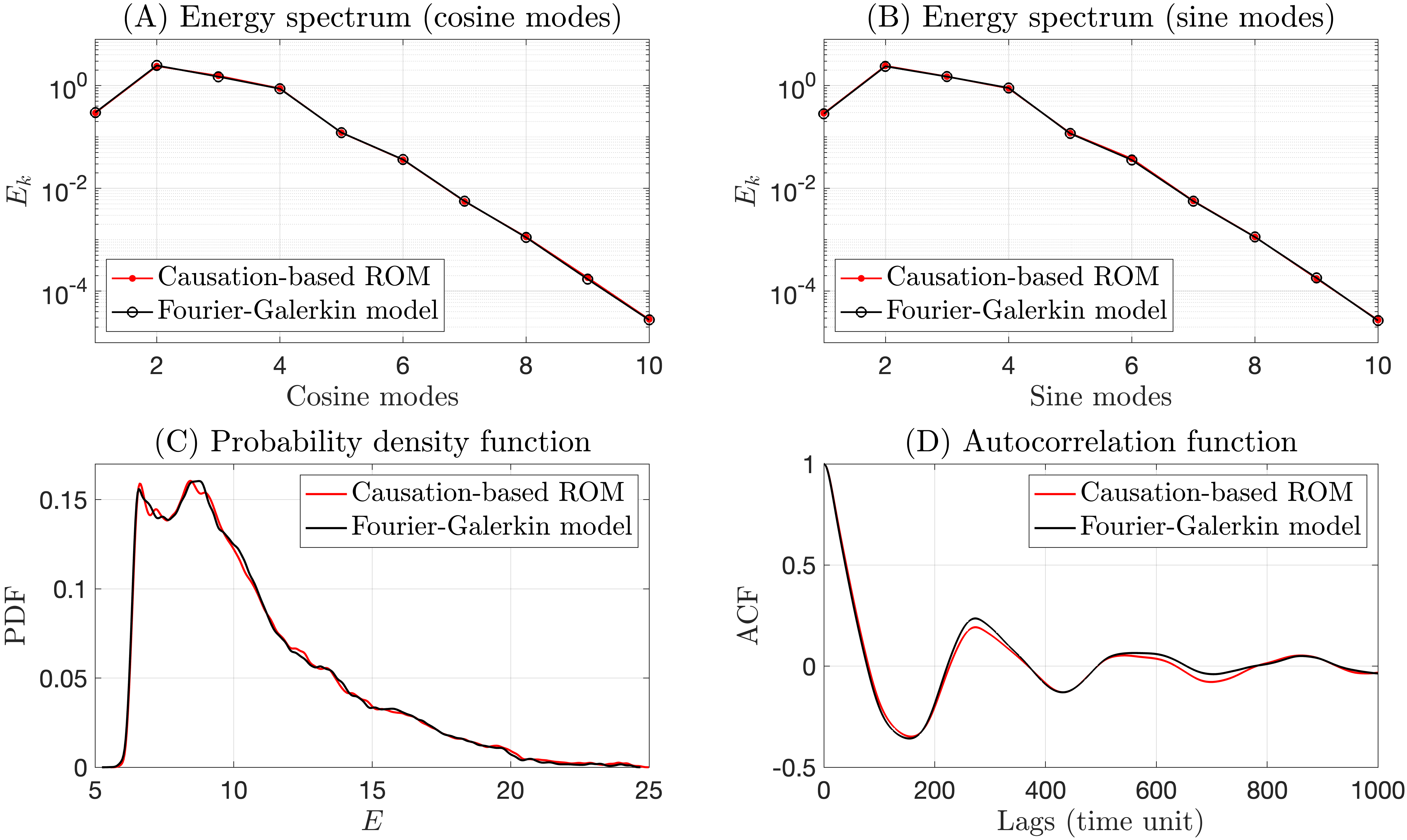}
\caption{\footnotesize Comparison of statistics between the 20-dimensional Fourier-Galerkin system and the associated causation-based ROM. The energy spectrum $E_k$ shown here is computed outside of the training window by averaging $(a^{\ell}_k(t))^2$ over the time window $[5\times 10^4, 3\times 10^5]$, for all the components $k=1,\ldots, 10$ and $\ell = 0, 1$. The cosine modes correspond to $\ell = 0$ (Panel A) and sine modes correspond to $\ell = 1$ (Panel B). The PDFs and the ACFs (shown in the bottom panels) for the kinetic energy $E(t) = \sum_{k=1}^{10} \sum_{\ell = 0}^1 (a_k^{\ell}(t))^2$ are computed over the same time window $[5\times 10^4, 3\times 10^5]$. We also note that the energy spectra for the sine and cosine modes corresponding to the same frequency $k$ are essentially the same for the Galerkin model (as well as for the learned model), indicating a type of equipartition of energy occurring in the model.}
\label{Fig_statistics_eigenbasis}
\end{figure}

With the constituent terms of the ROM identified now, we follow Step 4 of Sec.~\ref{Sec_CEM_overview} and use a standard MLE to determine the model coefficients; see Sec.~\ref{Subsec:ParameterEstimation}. For this purpose, we use the same training solution data used in the previous causal inference step. Among the 283 coefficients in $\bm{\Theta}$, 20 coefficients on the diagonal of the matrix are to be learned for the linear part, and the remaining 263 coefficients are for the nonlinear terms. The numerical values of these coefficients are graphically shown in Fig.~\ref{Fig_learned_coef_eigenbasis}. For the linear part (top panels of Fig.~\ref{Fig_learned_coef_eigenbasis}), the learned model coefficients recover these for the true model with high precision. We have checked that the relative error is below $0.06\%$ for all the 20 coefficients. For the nonlinear terms (bottom panels of Fig.~\ref{Fig_learned_coef_eigenbasis}), the largest differences between the learned model coefficients and the true ones occur at the 12 mismatched terms, as expected. Outside of these mismatched terms, the error is one-order smaller (on the scale of $10^{-3}$) compared with those shown in panel F of Fig.~\ref{Fig_learned_coef_eigenbasis}. Note also that the causation entropy for the mismatched terms all fall below the red dashed line in Fig.~\ref{Fig_CEM_coef_eigenbasis} and thus being filtered out in the learned model for this chosen cutoff threshold. We then expect that they play a less important role than the other 283 terms in ``orchestrating'' the dynamics in the true model. We also note that the learned noise amplitude matrix $\boldsymbol{\sigma}$ associated with the Gaussian noise term (see $\boldsymbol{\sigma} \dot{\bm{W}}(t)$ in \eqref{Identified_Model}) comes with very small entries on the scale of $10^{-7}$. Thus, the noise term is essentially negligible in the resulting causation-based ROMs. It turns out that the learned model with this cutoff threshold can already capture faithfully the true dynamics as shown in Figs.~\ref{Fig_soln_eigenbasis} and \ref{Fig_statistics_eigenbasis}.

In particular, we show in Fig.~\ref{Fig_soln_eigenbasis} the reconstructed spatiotemporal field for the true Galerkin model defined by \eqref{Eq_uG_eigen} (left panel) and its analog from the causation-based ROM (right panel). The solutions are shown in a time window that is far beyond the training window $[10^4, 5\times 10^4]$. The chaotic dynamics from the learned model are essentially indistinguishable from those in the true model, with local maxima (reddish patches) and local minima (bluish patches) progressing in a zigzag way as time evolves, forming a rich variety of local patterns. In that respect, we also point out that the long thin reddish strip observed in the left panel formed in the time window $[2.96\times 10^5, 2.97\times 10^5]$, which propagates from the left side of the domain ($x=0$) all the way up to almost the right side of the domain, has also been observed in other time windows for the learned model. This good reproduction of the dynamics is further confirmed at the statistical level as shown in Fig.~\ref{Fig_statistics_eigenbasis} for the energy spectrum $E_k$ (top panels) as well as the PDF and the autocorrelation function (ACF) of the kinetic energy $E$ (bottom panels); see the caption of this figure for further details.

Going back to Fig.~\ref{Fig_CEM_coef_eigenbasis}, when the more conservative cutoff threshold 0.05 (corresponding to the blue dashed line) is used, the corresponding causation-based ROM contains a total of 354 terms, which includes all of the 295 terms appearing in the true Galerkin model. The performance of this new causation-based ROM is similar to those shown in Fig.~\ref{Fig_soln_eigenbasis}B and Fig.~\ref{Fig_statistics_eigenbasis}. This indicates that the terms whose causation entropy values fall in between the gap marked by the red and blue dashed lines in Fig.~\ref{Fig_CEM_coef_eigenbasis} already play very little role in determining the dynamics of the learned model.

We also note that the existence of a clear gap to separate the larger and smaller causation entropy values, such as shown in Fig.~\ref{Fig_CEM_coef_eigenbasis} seems to be tied to the fact that the Fourier-Galerkin systems \eqref{KSE_Galerkin_eig} themselves admit a sparse structure; see again \eqref{Sparse1}--\eqref{Sparse3}. When other (global) basis functions are used, the corresponding Galerkin system may no longer be sparse. As such, one should no longer expect a clear gap to present in the causation entropy plot. However, as shown below using the POD basis, the ranking of the library terms provided by the causation entropy still offers a compelling way to obtain skillful yet significantly sparsified models.

\subsection{Data-driven inverse models under the POD basis} \label{Sec_POD_ROM_results}

We turn now to examine the situation when the underlying orthogonal basis is constructed empirically instead, which is taken to be the POD basis here. For benchmarking purposes, we will compare the performance of the learned model with that of the POD-Galerkin system \eqref{KSE_Galerkin_POD} as well as a thresholded version of the Galerkin system obtained by removing terms whose coefficients in absolute value are below a given threshold to achieve a specified sparsity percentage.

\begin{figure}
\centering
\includegraphics[width=1\textwidth]{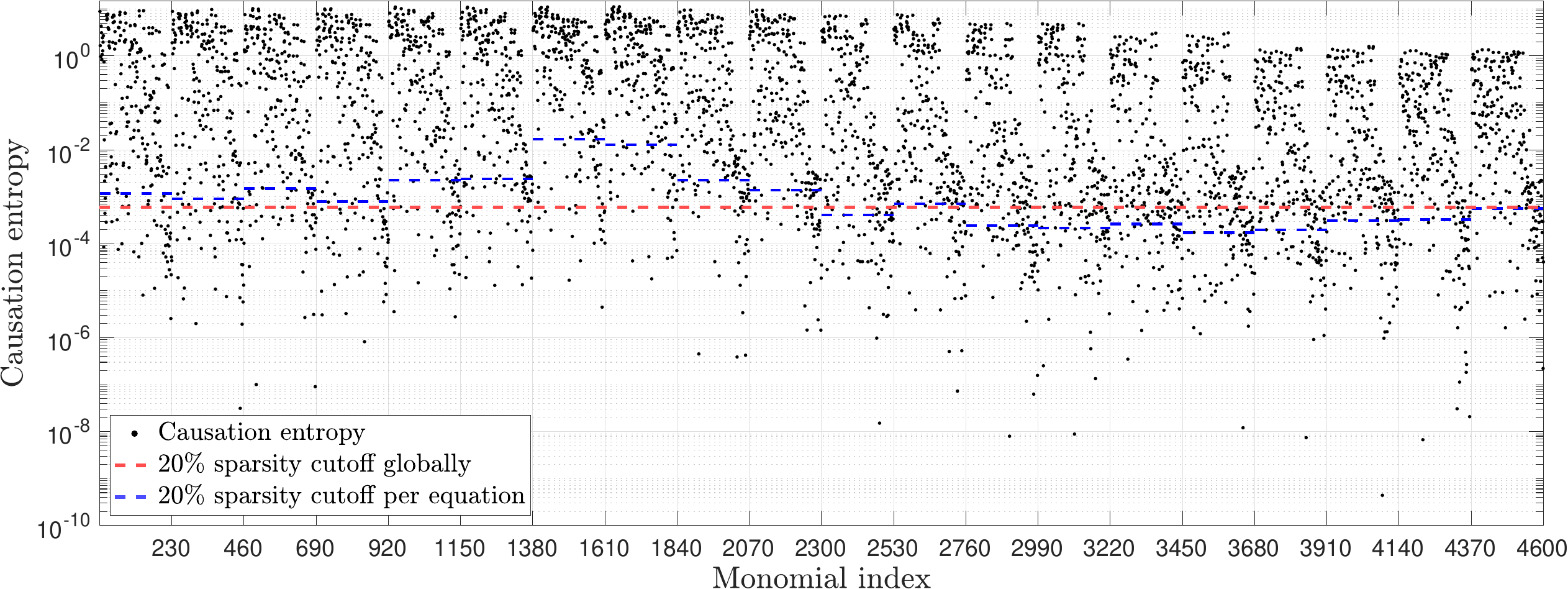}
\caption{\footnotesize Causation entropy that ranks the library functions for learning a data-driven quadratic inverse model of the 20-dimensional POD Galerkin system \eqref{KSE_Galerkin_POD}. The causation entropy values are grouped by equation in the same way as done in Fig.~\ref{Fig_CEM_coef_eigenbasis}. Also shown are the cutoff thresholds for ensuring 20\% sparsity based on two strategies: the red dashed line corresponds to the threshold $6\times 10^{-4}$ that separates the lower 20\% of all the 4600 causation entropy values from the remaining 80\%, while the blue dashed line segments correspond to the thresholds that separate the lower 20\% of the 230 causation entropy values for each of the 20 equations. The total number of terms kept in the learned models based on these two strategies are the same, but the constituent terms kept in the corresponding identified models are slightly different.}
\label{Fig_CEM_coef_POD}
\end{figure}

\begin{figure} 
\centering
\includegraphics[width=1\textwidth]{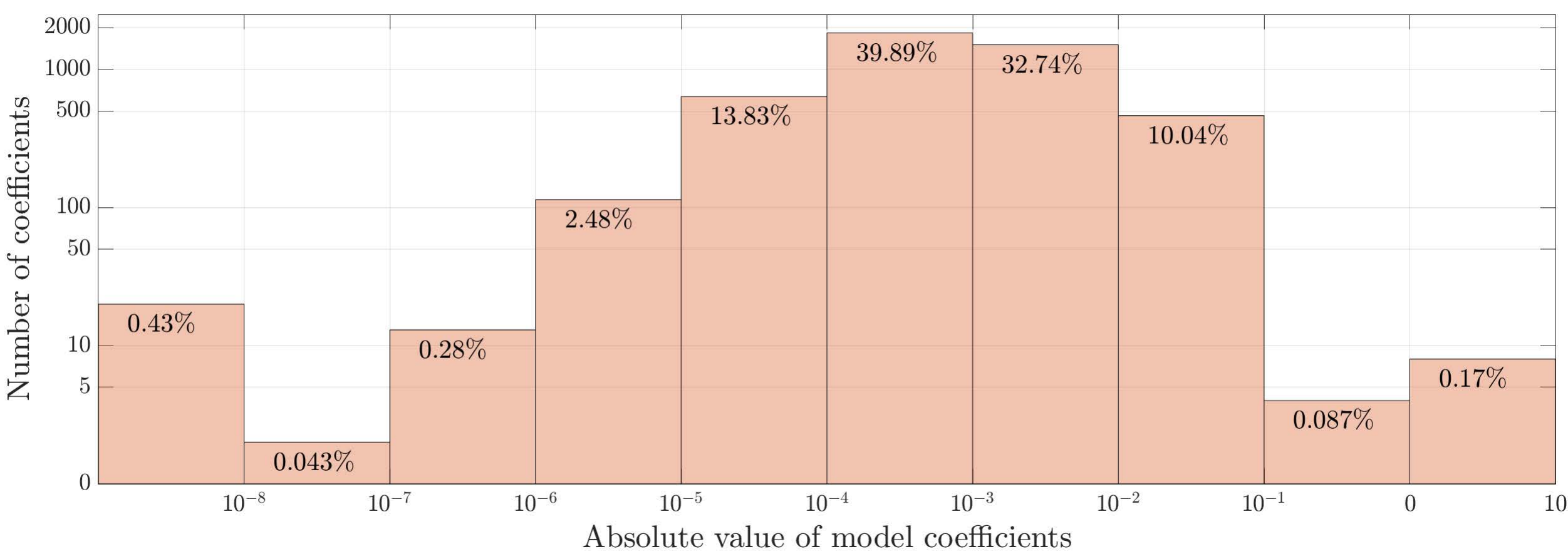}
\caption{\footnotesize Distribution of the model coefficients for all the 4600 terms in the 20-dimensional POD-Galerkin system \eqref{KSE_Galerkin_POD} of the KSE \eqref{Eq_KSE}, for the parameter regime given by Sec.~\ref{Sec_numerical_setup}.}
\label{Fig_POD_Galerkin_coef_dist}
\end{figure}

The causation entropy, as computed using the 20-dimensional POD projection of the KSE solution, is shown in Fig.~\ref{Fig_CEM_coef_POD}. Unlike the case with the Fourier basis shown in Fig.~\ref{Fig_CEM_coef_eigenbasis}, we no longer see a gap that separates a small fraction of larger causation entropy values with the remaining smaller causation entropy values. As mentioned at the end of the previous subsection, a plausible reason is that the POD-Galerkin system \eqref{KSE_Galerkin_POD} itself does not have a sparse structure. Recall that the 20-dimensional Fourier-Galerkin system \eqref{KSE_Galerkin_eig} utilized in the previous subsection has only 295 terms in its vector field, accounting for about $6.41\%$ of the total $20\times 230 = 4600$ possible monomials in a 20-dimensional quadratic vector field (excluding constant terms). In sharp contrast, almost all the 4600 terms are present in the 20-dimensional POD-Galerkin system \eqref{KSE_Galerkin_POD}. As shown in Fig.~\ref{Fig_POD_Galerkin_coef_dist}, the absolute value of the coefficients falls in the range $[10^{-5}, 10^{-1}]$ for $96.5\%$ of the terms (namely 4439 terms) in this POD-Galerkin system.

Due to the lack of any obvious cutoff thresholds appearing in the distribution of the causation entropy values, a possible way to proceed is to construct a hierarchy of inverse models that maintain different sparsity percentages by adjusting the cutoff threshold. There are two different ways to carry out this cutoff procedure. One way is to choose a uniform cutoff threshold for all the equations, such as indicated by the red dashed line in Fig.~\ref{Fig_CEM_coef_POD}. Apparently, this approach only ensures that a given percentage of terms is removed from the learned model but does not guarantee that the percentage of terms removed is the same for each equation in the system. The other way is to choose a custom cutoff threshold for each equation, such as indicated by the blue dashed lines in Fig.~\ref{Fig_CEM_coef_POD} to achieve the same sparsity percentage for each equation. In principle, the two cutoff procedures can lead to quite different reduced models, especially when the range of the causation entropy values varies significantly from equation to equation. However, for the model considered here, it has been checked that the ROMs obtained by the two approaches for a given sparsity percentage lead to similar modeling performance. For all the numerical results reported below, the causation-based ROMs are constructed using the latter approach to gain the same sparsity percentage for all the equations. Once the constituent terms are determined based on a chosen cutoff threshold strategy for the causation entropy values, we use again the MLE to determine the model coefficients in the causation-based ROMs; see Sec.~\ref{Subsec:ParameterEstimation}.

\begin{figure}[h]
\centering
\includegraphics[width=1\textwidth]{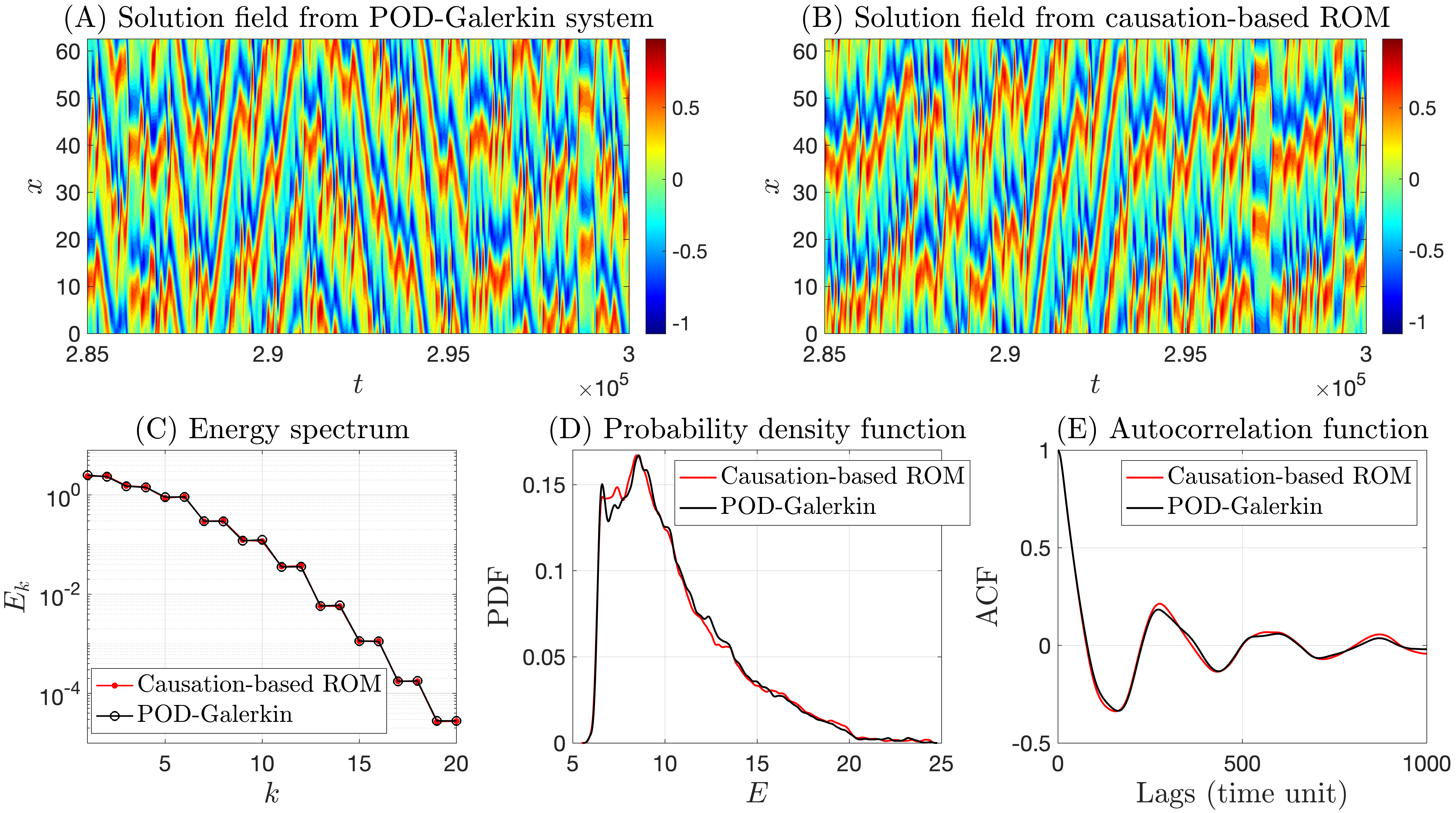}
\caption{\footnotesize Performance of the 20-dimensional causation-based ROM under POD basis, with a 20\% sparsity cutoff per equation, in comparison with the 20-dimensional POD-Galerkin system. The reconstructed solution fields are shown in Panel A for the POD-Galerkin system and Panel B for the learned model. The energy spectrum $E_k$ is shown in Panel C. The PDFs and the ACFs of the kinetic energy $E$ are shown in Panels D and E, respectively. The energy spectra $E_k$'s and the kinetic energy $E$ are computed in the same way as described in the caption of Fig.~\ref{Fig_statistics_eigenbasis}. Like in Fig.~\ref{Fig_statistics_eigenbasis}, the time window used for computing $E_k$ and $E$ is $[5\times 10^4, 3\times 10^5]$, which is outside of the training window $[10^4, 5\times 10^4]$.}
\label{Fig_soln_PODbasis}
\end{figure}

In Fig.~\ref{Fig_soln_PODbasis}, we present the skill of the 20-dimensional causation-based ROM with a 20\% sparsity. As can be observed, even though the ROM contains 20\% fewer terms than the corresponding POD-Galerkin system, it can faithfully reproduce the essential dynamical features and the associated statistics appearing in the solution of the POD-Galerkin system. We also checked that even by increasing the sparsity percentage to 50\%, the causation-based ROM can still produce chaotic transient dynamics over a long time window with the corresponding solution field resembling that shown in Panels A and B of Fig.~\ref{Fig_soln_PODbasis}, although the solution eventually becomes periodic after about $3.6\times 10^6$ time-step iterations. This slow drift to periodic dynamics after long-time integration observed for ``severely'' truncated causation-based ROMs is not a surprise. The KSE is known to have many periodic dynamics regimes interlaced with chaotic dynamics regimes \cite{hyman1986order}. In other words, the chaotic attractors observed for the KSE may be prone to instability under perturbations depending on the parameter regimes considered. Since a ROM can be viewed as a perturbation of the original KSE model, it is possible for the dynamics of a highly truncated ROM to be (gradually) pushed towards the basin of attraction of a periodic attractor in a nearby regime.

Although the dynamical features of the KSE dictate that one may not be able to use a too sparse ROM to capture long-term statistics for certain parameter regimes, the fact that such a causation-based ROM can still reproduce accurately short-time features suggests its potential usage for other purposes such as data assimilation and short-term trajectory prediction. In the next subsection, we demonstrate the advantage of such highly sparse causation-based ROMs in the context of data assimilation with partial observations.

\subsection{Application to data assimilation with partial observations} \label{Sec_DA_results}

We now illustrate the performance of causation-based ROMs in the context of data assimilation to recover unobserved higher-frequency mode dynamics based on observation data for a few low-frequency mode dynamics. As a benchmark, we also compare the results obtained from a thresholded stochastic POD-Galerkin model described below. For simplicity, we assume that the observation data is available continuously in time, and we perform the data assimilation with the ensemble Kalman-Bucy filter (EnKBF) \cite{bergemann2012ensemble, amezcua2014ensemble} for both reduced systems.

The thresholded Galerkin system is obtained from the true 20-dimensional POD-Galerkin system as follows. We rank the coefficients of the true Galerkin system from large to small in absolute value and then drop the terms with coefficients below a cutoff value that is determined to ensure the number of monomials retained is the same as that of the employed causation-based ROM. After identifying the terms to be kept, we use the MLE to estimate the model coefficients and the covariance matrix of the additive noise term in the final thresholded Galerkin model.

For the sake of clarity, we first provide below some details about the EnKBF applied to a generic $n$-dimensional SDE system of the form \eqref{Identified_Model}:
\be \label{Identified_Model_recall}
 \frac{\d \bm{a}}{\d t} = \boldsymbol{\Phi}(\bm{a}) + \boldsymbol{\sigma} \dot{\bm{W}}(t),
\ee
 in which the noise amplitude matrix $\boldsymbol{\sigma}$ is assumed to be $n \times n$-dimensional and the first $r$ component of $\bm{a}$ is taken to be observed while the remaining components to be unobserved. We denote
\bea
& \bm{y} = (a_1, \ldots, a_r)^{\mathtt{T}}, && \bm{z} = (a_{r+1}, \ldots, a_n)^{\mathtt{T}}, \\
& \bm{W}_1 = (W_1, \ldots, W_r)^{\mathtt{T}}, && \bm{W}_2 = (W_{r+1}, \ldots, W_n)^{\mathtt{T}}.
\eea
We also decompose $\boldsymbol{\sigma}$ into four submatrices
\be
\boldsymbol{\sigma} =
\begin{pmatrix}
\boldsymbol{\sigma}_{11} & \boldsymbol{\sigma}_{12} \\
\boldsymbol{\sigma}_{21} & \boldsymbol{\sigma}_{22}
\end{pmatrix},
\ee
where the dimensions of $\boldsymbol{\sigma}_{11}$ and $\boldsymbol{\sigma}_{22}$ are $r\times r$ and $(n-r)\times (n-r)$, respectively. We then rewrite \eqref{Identified_Model_recall} using $(\bm{y}, \bm{z})$ as follows
\bea \label{SDE_rewritten}
 & \frac{\d \bm{y}}{\d t} = \bm{g}_1(\bm{y}, \bm{z}) + \boldsymbol{\sigma}_{11} \dot{\bm{W}}_1(t),  \\
 & \frac{\d \bm{z}}{\d t} = \bm{g}_2(\bm{y}, \bm{z}) + \boldsymbol{\sigma}_{22} \dot{\bm{W}}_2(t),
\eea
where $\bm{g}_1$ and $\bm{g}_2$ denote respectively the first $r$ and the remaining $n-r$ components of the (nonlinear) function $\Phi$ in \eqref{Identified_Model_recall}. Note also that compared with the original system \eqref{Identified_Model_recall}, we decoupled the noise terms in the $\bm{y}$- and $\bm{z}$-subsystems by dropping $\boldsymbol{\sigma}_{12} \dot{\bm{W}}_2(t)$ in the $\bm{y}$-subsystem and $\boldsymbol{\sigma}_{21} \dot{\bm{W}}_1(t)$ in the $\bm{z}$-subsystem for simplicity. 
In practice, the noise amplitude matrix $\boldsymbol{\sigma}$ is oftentimes diagonally dominant. This is, in particular, true for the KSE problem considered here. Additionally, the noise in both the causation-based ROM and the thresholded POD-Galerkin system is very weak. Thus, such an approximation has little impact on the accuracy of final data assimilation results.

Assume that a total of $p$ ensemble members are used in the EnKBF. Denote the collection of all the $p$ ensemble members at time $t$ by
\bes
\bm{Z}(t) = (\bm{z}_1(t), \bm{z}_2(t), \ldots, \bm{z}_p(t))^{\mathtt{T}}.
\ees
Denote also the observation data of $\bm{y}$ at time $t$ by $\bm{y}_{\text{obs}}(t)$. We define then
\begin{equation}
\overline{\bm{Z}}(t) = \frac{1}{p}\sum_{\ell=1}^p  \bm{z}_\ell(t), \quad  \overline{\bm{g}}_2 (\bm{y}_{\text{obs}}(t), \bm{Z}(t)) = \frac{1}{p} \sum_{\ell=1}^p \bm{g}_2(\bm{y}_{\text{obs}}(t), \bm{z}_{\ell}(t)),
\end{equation}
and
\begin{equation}\label{EnKBF_Variance}
\begin{split}
& \mathcal{N}(\bm{y}_{\text{obs}}(t),\bm{Z}(t)) = \frac{1}{(p-1)} \sum_{\ell=1}^p (\bm{z}_{\ell}(t) - \overline{\bm{Z}}(t))\Big( \bm{g}_2(\bm{y}_{\text{obs}}(t), \bm{z}_{\ell}(t)) - \overline{\bm{g}}_2(\bm{y}_{\text{obs}}(t), \bm{Z}(t)) \Big)^{\mathtt{T}} C^{-1},
\end{split}
\end{equation}
where $C = \boldsymbol\sigma_{22}\boldsymbol\sigma_{22}^\mathtt{T}$.

Then, each ensemble member $\bm{z}_i$, $i=1,2,\ldots,p$, of the EnKBF is computed using
\begin{equation} \label{Eq_EnKBF}
\begin{split}
\frac{\d \bm{z}_i}{\d t} & = \bm{g}_2(\bm{y}_{\text{obs}}(t),\bm{z}_i) + \boldsymbol{\sigma}_{22} \dot{\bm{W}}_{2,i}(t)\\
&\qquad - \mathcal{N}(\bm{y}_{\text{obs}}(t),\bm{Z}(t))\big[\bm{g}_1(\bm{y}_{\text{obs}}(t),\bm{z}_i) - \dot{\bm{y}}_{\text{obs}}(t) + \boldsymbol{\sigma}_{11} \dot{\bm{W}}_{1,i}(t) \big],
\end{split}
\end{equation}
where $\bm{W}_{1,i}$ and $\bm{W}_{2,i}$ are respectively $r$-dimensional and $(n-r)$-dimensional Brownian motions for $i=1,2,\ldots,p$, with their
components to be all mutually independent.

The setup of the data assimilation experiment is as follows. We observe the amplitudes of the first three POD modes of the KSE solutions and aim to recover the amplitudes of the few dominant unobserved modes by applying the EnKBF to either the 20-dimensional causation-based ROM with a large sparsity percentage or the corresponding 20-dimensional thresholded POD-Galerkin system with the same sparsity percentage. We take the size of the EnKBF ensemble to be $p=500$, and the unobserved variables are initialized to be zero for all the ensemble simulations. The sparsity percentage of the ROMs is taken to be 90\%, resulting in 460 terms in the drift part of both the causation-based ROM and the thresholded POD-Galerkin system. The KSE is simulated over the time window $[0, 2000]$ with the initial data taken to be the solution profile at the last time instant of the training data utilized for constructing the POD basis function as well as the training of the ROMs.

\begin{figure}
\centering
\includegraphics[width=1\textwidth]{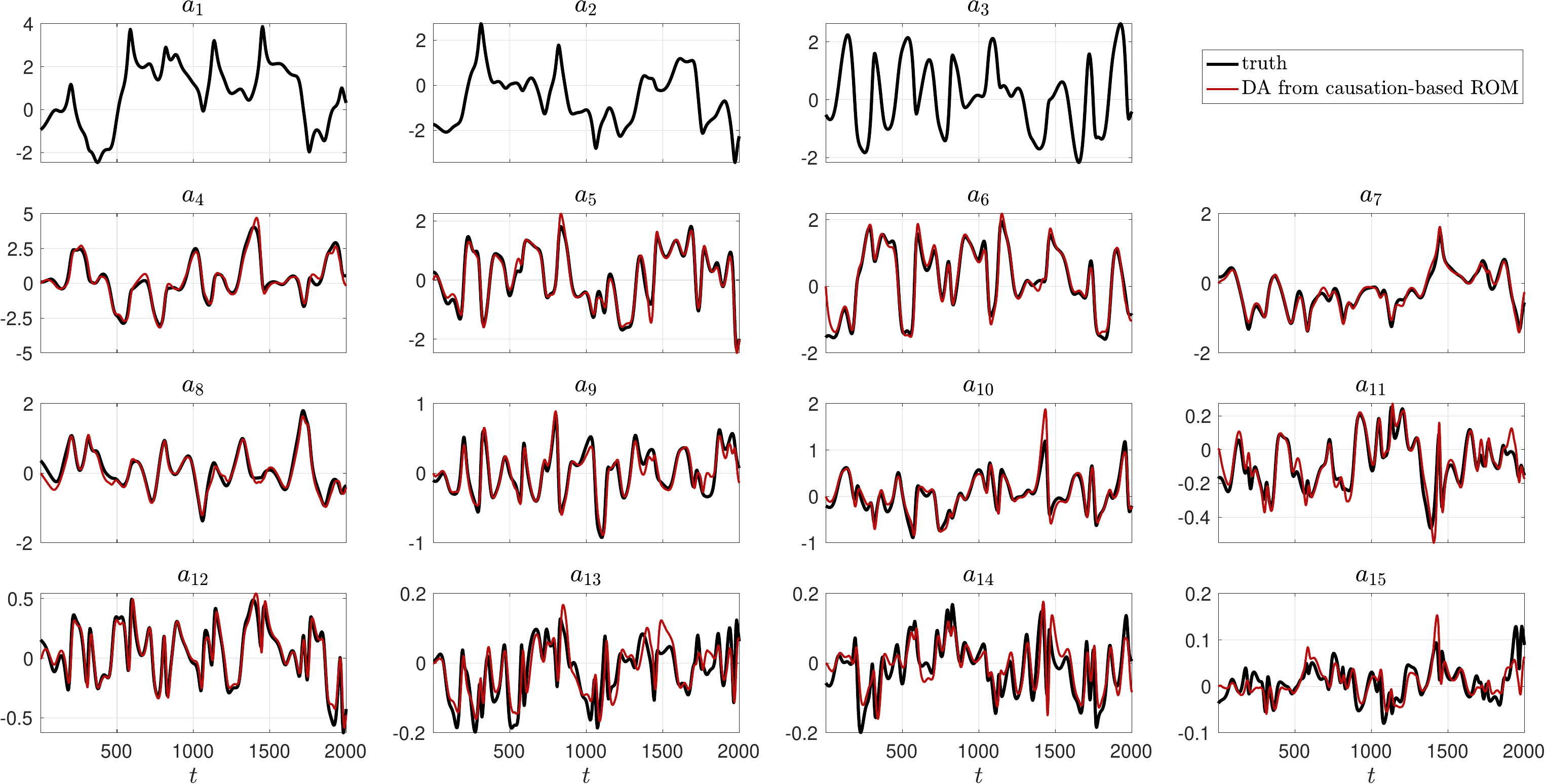}
\caption{\footnotesize The black curves show the projections of the true KSE solution onto the first 15 POD modes, with the first three components $a_1$, $a_2$, and $a_3$ taken to be the observed modes in the data assimilation experiments. The red curves show the assimilated ensemble mean dynamics of the unobserved POD modes from the 20-dimensional causation-based ROM, with 90\% sparsity.} \label{Fig_DA_cem_timeseries}
\end{figure}

In the first row of Fig.~\ref{Fig_DA_cem_timeseries}, we show the time series of the three observed POD modes. On average, these three modes capture about 63.5\% of the kinetic energy in the KSE solution for the considered parameter regime, while above 99\% of the kinetic energy is captured by the first 10 POD modes. As shown in Fig.~\ref{Fig_DA_cem_timeseries} (black curves), modes 4 to 10 still have quite large amplitude oscillations almost comparable with those of the first three modes, and they evolve on different time scales. The fact that the unobserved dynamics still contain, on average, nearly 40\% of the kinetic energy and that their projected dynamics reveal multi-scale, highly chaotic oscillatory features present arguably a challenging test ground for the data assimilation experiment.

The red curves in Fig.~\ref{Fig_DA_cem_timeseries} represent the posterior mean of the unobserved modes $a_4$ through $a_{15}$ (red curves) obtained from the EnKBF applied to the causation-based ROM, in comparison with the corresponding true POD projections of the KSE solutions (black curves). Despite its highly sparse nature, with 90\% sparsity compared with the true POD-Galerkin system of the same dimension, the causation-based ROM is able to recover with high fidelity all the energetic unobserved modes, $a_4, \ldots, a_{12}$. The skill for the remaining small amplitude modes, $a_{13}, \ldots, a_{20}$, deteriorates as the mode index increases, as can be seen in Fig.~\ref{Fig_DA_cem_timeseries} for modes $a_{13}$, $a_{14}$, and $a_{15}$. However, these 8 modes contain, on average, only approximately $0.13\%$ of the solution energy.

\begin{figure}
\centering
\includegraphics[width=1\textwidth]{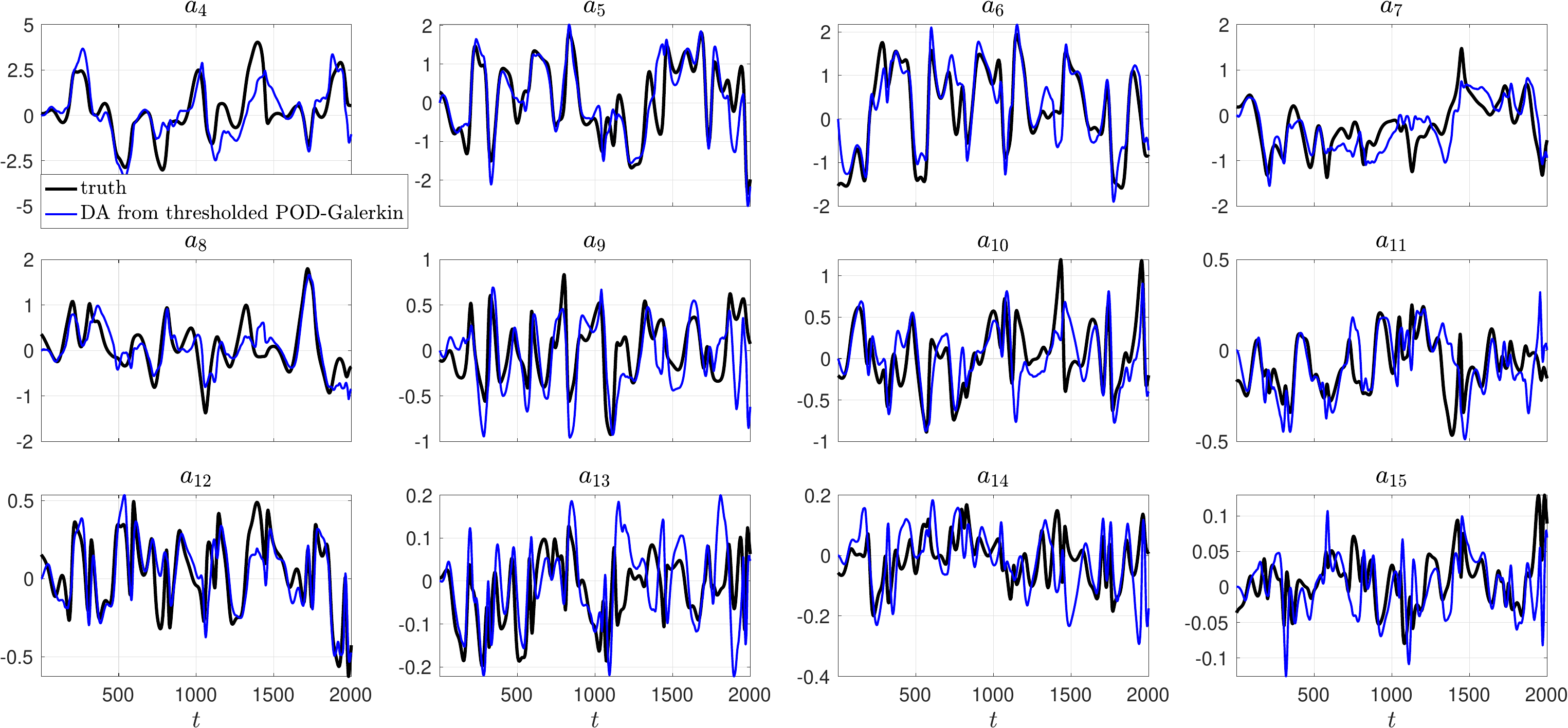}
\caption{\footnotesize The assimilated ensemble mean dynamics of the unobserved POD modes from the EnKBF of the 20-dimensional thresholded POD-Galerkin system with 90\% sparsity (blue curve), in comparison with the POD projections of the true KSE solution (black curves).} \label{Fig_DA_thresholded_Galerkin_timeseries}
\end{figure}

\begin{figure}[tbh!]
\centering
\includegraphics[width=0.8\textwidth]{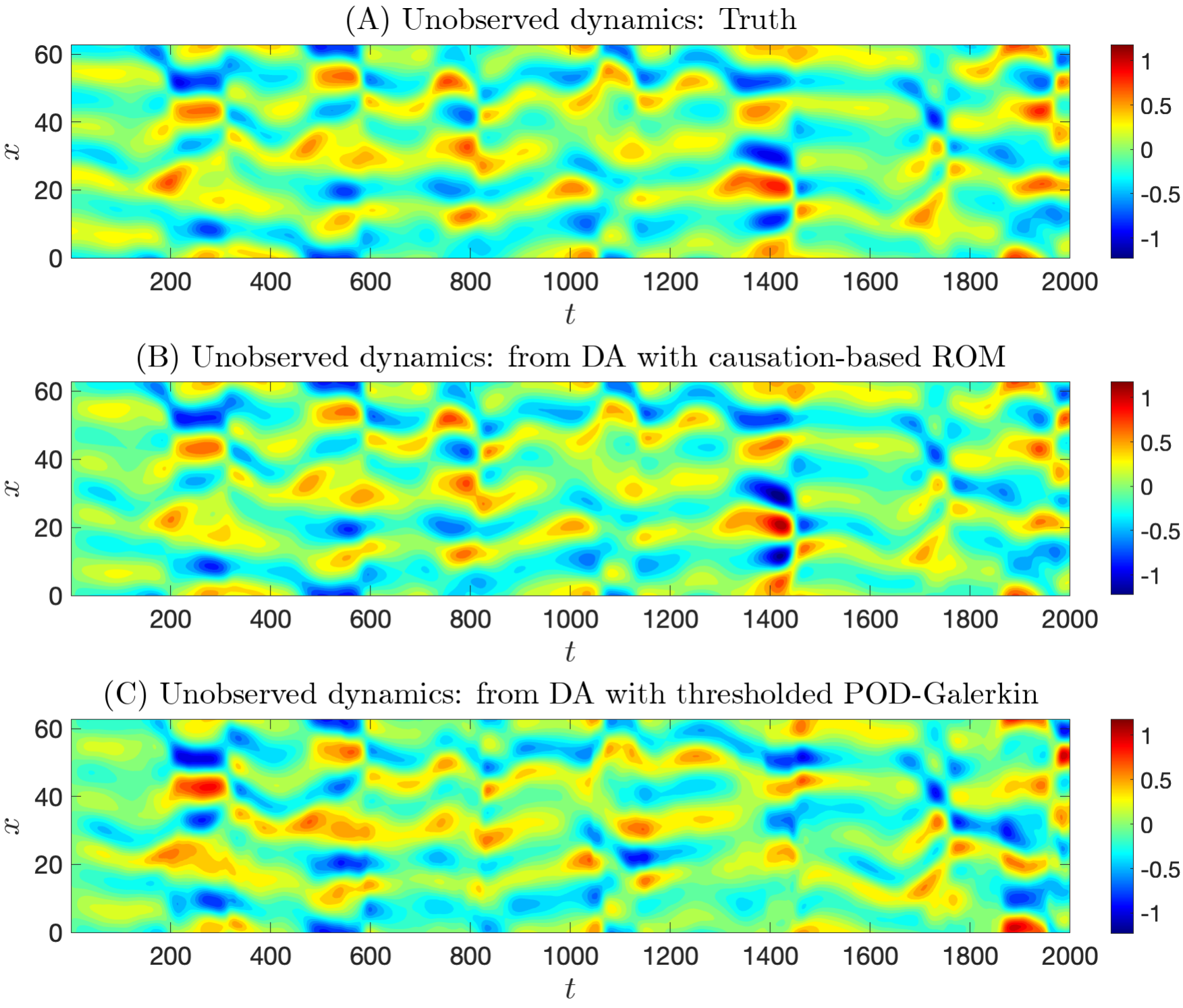}
\caption{\footnotesize Comparison of the true unobserved spatiotemporal field (top panel) with those reconstructed based on the assimilated ensemble mean dynamics shown in Fig.~\ref{Fig_DA_cem_timeseries} for the causation-based ROM (middle panel) and in Fig.~\ref{Fig_DA_thresholded_Galerkin_timeseries} for the thresholded POD-Galerkin system (bottom panel).} \label{Fig_DA_fields}
\end{figure}

As a comparison, the corresponding results for the 20-dimensional thresholded POD-Galerkin system with 90\% sparsity are shown in Fig.~\ref{Fig_DA_thresholded_Galerkin_timeseries}. The skill is significantly worse than that obtained from the causation-based ROM. When comparing these unobserved dynamics at the spatiotemporal field level, it also reveals that the thresholded POD-Galerkin system with this high truncation ratio suffers particularly severely when there is a relatively abrupt change in the solution dynamics, such as shown at around $t = 1450$ in Fig.~\ref{Fig_DA_fields}. Finally, we mention that the time series for all the 500 ensemble members in the data assimilation essentially coincide with each other for both of the two ROMs analyzed due to the fact that the involved noise amplitude matrix $\boldsymbol{\sigma}$ for both ROMs employed has entries all close to zero.

The above results show that causation entropy can indeed be utilized to rank the relative importance of candidate terms from a given function library for the construction of skillful sparse inverse models. The obtained superior data assimilation skills compared with those from the thresholded POD-Galerkin system also illustrate that a naive truncation based on the numerical values of the model coefficients in e.g.~a POD-Galerkin system may not be appropriate, especially when a highly truncated ROM is sought.

\section{Discussion and Conclusions}\label{Conclusion}

In this article, we analyzed an efficient approach to identifying data-driven ROMs with a sparse structure using a quantitative indicator called causation entropy. For each potential building-block function $f$ in the vector field of the $i$-th component $a_i$, the associated causation entropy measures the difference between the entropy of $\dot{a}_i$ conditioned on the whole set of candidate functions and the one conditioned on the set without $f$; see \eqref{Causation_Entropy}. Thus, it quantifies statistically the additional contribution of each term to the underlying dynamics beyond the information already captured by all the other terms in the model ansatz.

The ranking of the candidate terms provided by the causation entropy leads to a hierarchy of parsimonious structures for the ROMs controlled by a cutoff threshold parameter. The model coefficients for the corresponding causation-based ROMs can then be learned using standard parameter estimation techniques, such as the MLE; cf.~Sec.~\ref{Subsec:ParameterEstimation}. Illustrating on the Kuramoto-Sivashinky equation, we showed in Sec.~\ref{Sec:KSE} that the obtained causation-based ROMs are skillful in both recovering long-term statistics and inferring unobserved dynamics via data assimilation when only a small subset of the ROM's state variables is observed.

We conclude by outlining some potential future directions to be explored. For this purpose, we want to emphasize first that, when building up the causation-based ROMs, it is straightforward to add additional physically relevant constraints, such as skew symmetry for certain linear terms and energy conservation for the quadratic nonlinearity. For the results shown in Sec.~\ref{Subsec:ParameterEstimation}, the obtained ROMs turn out to be stable without enforcing energy conservation constraints, even though the quadratic term in the Kuramoto-Sivashinky equation conserves energy. However, such a constraint is expected to be important, e.g., in the reduction of fluid problems in turbulent regimes. To enforce such constraints, we just need to make sure all relevant terms are included in the identified model structure since it can happen that the causation entropy for some but not all of the terms involved in the constraint is above a given cutoff threshold.
Of course, the subsequent parameter estimation is subject to the desired constraints as well, which can be performed using, e.g., the constrained MLE \cite[Section 2.5]{chen2023causality}.

Oftentimes, when constructing ROMs for highly chaotic systems, one needs to include closure terms to take into account the impact of the orthogonal dynamics not resolved by the ROMs \cite{ahmed2021closures}. Different strategies can be envisioned to extend the current framework for this purpose. For instance, after the drift part of the causation-based ROM is identified (i.e., the $\boldsymbol\Phi(\bm{a})$-term in \eqref{Identified_Model} or \eqref{RHS_Model}), instead of fitting the resulting training residual data by a noise term $\boldsymbol{\sigma} \dot{\bm{W}}(t)$, we can explore more advanced data-driven techniques such as multilevel approaches and empirical model reduction \cite{kondrashov2015data,kravtsov2005multilevel,majda2012physics}, nonlinear autoregressive techniques \cite{lu2017data,billings2013nonlinear,chen2022shock,chorin2015discrete}, or neural networks. Alternatively, one could first learn a higher-dimensional causation-based ROM, then use parameterization techniques \cite{CLM20,chekroun2021stochastic} to approximate the newly added components by those to be kept.

To what extent one can sparsify a ROM depends apparently on the purposes of the ROMs. However, it can also be tied to the underlying orthogonal basis employed. As already seen in Sec.~\ref{Sec:KSE}, the causation-based ROMs constructed using the eigenbasis come with a much sparser structure than those built from the POD basis, for the PDE considered. It would be interesting to explore if a coordinate transformation exists that can further enhance the sparsity of the ROMs built on a POD basis. For instance, if we rewrite the POD-ROM under the eigenbasis of the ROM's linear part, we can oftentimes achieve a diagonalization of the linear terms since eigenvalues with multiplicity one are generic. However, whether this transformation can also help aggregate the nonlinearity to form sparser structures (after re-computing the causation entropy matrix under the transformed basis) is up to further investigation.

Another aspect concerns the efficient computation of the causation entropy matrix when the number of functions, $M$, in the learning library is in the order of several thousand or beyond, which can, for instance, be encountered for ROMs with dimension $100$ or higher. The computational cost for determining
a causation entropy lies in the calculation of the log-determinants of the four covariance matrices involved in formula \eqref{Entropy_Gaussians}, which are of dimension either $M\times M$ or $(M \pm 1)\times (M \pm 1)$. For a ROM of dimension $N$, there are a total of $N\times M$ causation entropies to determine. Thus, we need to compute the log-determinants of $4\times N \times M$ covariance matrices, each with dimension about $M\times M$.
To gain computational efficiency when $M$ is large, one may benefit from techniques for approximating the log-determinant of a high dimensional symmetric positive definite matrix \cite{pace2009sampling,boutsidis2017randomized}, although additional investigation would be needed to see how one can strike a balance between the computational efficiency gained and the approximation error made on each entry of the causation entropy matrix. Alternatively, we can try to reduce the number of functions in the feature library either by an iterative approach using a greedy algorithm \cite{sun2015causal} or by exploring potential physical/modeling insights for the considered applications. For instance, in \cite{chen2023causality}, a localization strategy is introduced to significantly reduce the size of the feature library when constructing an efficient causation-based ROM for the two-layer Lorenz 1996 model.

\section*{Acknowledgments}
This research was funded in part by the Army Research Office grant W911NF-23-1-0118 (N.C.), the Office of Naval Research grant N00014-24-1-2244 (N.C.),  and the National Science Foundation grants DMS-2108856 and DMS-2407483 (H.L.). We also acknowledge the computational resources provided by Advanced
Research Computing at Virginia Tech. A preprint of this manuscript is available on arxiv.org \cite{chen2024minimum}. 

\section*{Conflict of Interest}
The authors have no conflicts to disclose.

\section*{Data Availability}
The data that support the findings of this study are available from the corresponding author upon reasonable request.


\begin{thebibliography}{100}

\bibitem{ahmed1989entropy}
N.~A. Ahmed and D. Gokhale.
\newblock Entropy expressions and their estimators for multivariate
  distributions.
\newblock {\em IEEE Transactions on Information Theory}, 35:688--692, 1989.

\bibitem{ahmed2021closures}
S.~E. Ahmed, S. Pawar, O. San, A. Rasheed, T. Iliescu, and B.~R. Noack.
\newblock {On closures for reduced order models--A spectrum of first-principle
  to machine-learned avenues}.
\newblock {\em Physics of Fluids}, 33(9):091301, 2021.

\bibitem{akaike1973information}
H. Akaike.
\newblock Information theory and an extension of the maximum likelihood
  principle.
\newblock In {\em Proceeding of the Second International Symposium on
  Information Theory}, pages 267--281. Akademiai Kiado, Budapest, 1973.

\bibitem{akaike1974new}
H. Akaike.
\newblock A new look at the statistical model identification.
\newblock {\em IEEE Transactions on Automatic Control}, AC-19:716--723, 1974.

\bibitem{almomani2020erfit}
A.~A. AlMomani and E. Bollt.
\newblock {ERFit}: Entropic regression fit {MATLAB} package, for data-driven
  system identification of underlying dynamic equations.
\newblock {\em arXiv preprint arXiv:2010.02411}, 2020.

\bibitem{almomani2020entropic}
A.~A.~R. AlMomani, J. Sun, and E. Bollt.
\newblock How entropic regression beats the outliers problem in nonlinear
  system identification.
\newblock {\em Chaos: An Interdisciplinary Journal of Nonlinear Science},
  30(1), 2020.

\bibitem{amezcua2014ensemble}
J. Amezcua, K. Ide, E. Kalnay, and S. Reich.
\newblock Ensemble transform {K}alman--{B}ucy filters.
\newblock {\em Quarterly Journal of the Royal Meteorological Society},
  140(680):995--1004, 2014.

\bibitem{armbruster1992phase}
D. Armbruster, R. Heiland, E.~J. Kostelich, and B. Nicolaenko.
\newblock {Phase-space analysis of bursting behavior in Kolmogorov flow}.
\newblock {\em Physica D: Nonlinear Phenomena}, 58:392--401, 1992.

\bibitem{aubry1993preserving}
N. Aubry, W.-Y. Lian, and E.~S. Titi.
\newblock Preserving symmetries in the proper orthogonal decomposition.
\newblock {\em SIAM Journal on Scientific Computing}, 14:483--505, 1993.

\bibitem{bergemann2012ensemble}
K. Bergemann and S. Reich.
\newblock {An ensemble {K}alman-Bucy filter for continuous data assimilation}.
\newblock {\em Meteorologische Zeitschrift}, 21:213--219, 2012.

\bibitem{bertozzi1998long}
A.~L. Bertozzi and M.~C. Pugh.
\newblock Long-wave instabilities and saturation in thin film equations.
\newblock {\em Communications on Pure and Applied Mathematics}, 51:625--661,
  1998.

\bibitem{bhola2023estimating}
S. Bhola and K. Duraisamy.
\newblock {Estimating global identifiability using conditional mutual
  information in a Bayesian framework}.
\newblock {\em Scientific Reports}, 13:18336, 2023.

\bibitem{billings2013nonlinear}
S.~A. Billings.
\newblock {\em {Nonlinear System Identification: {NARMAX} Methods in the Time,
  Frequency, and Spatio-temporal Domains}}.
\newblock John Wiley \& Sons, 2013.

\bibitem{Boers_al17}
N. Boers, M.~D. Chekroun, H. Liu, D. Kondrashov, D.-D. Rousseau, A. Svensson,
  M. Bigler, and M. Ghil.
\newblock {Inverse stochastic-dynamic models for high-resolution Greenland
  ice-core records}.
\newblock {\em Earth System Dynamics}, 8:1171--1190, 2017.

\bibitem{boninsegna2018sparse}
L. Boninsegna, F. N{\"u}ske, and C. Clementi.
\newblock Sparse learning of stochastic dynamical equations.
\newblock {\em The Journal of chemical physics}, 148(24), 2018.

\bibitem{boutsidis2017randomized}
C. Boutsidis, P. Drineas, P. Kambadur, E.-M. Kontopoulou, and A. Zouzias.
\newblock A randomized algorithm for approximating the log determinant of a
  symmetric positive definite matrix.
\newblock {\em Linear Algebra and its Applications}, 533:95--117, 2017.

\bibitem{branicki2012quantifying}
M. Branicki and A.~J. Majda.
\newblock Quantifying uncertainty for predictions with model error in
  non-{G}aussian systems with intermittency.
\newblock {\em Nonlinearity}, 25(9):2543, 2012.

\bibitem{brunton2016discovering}
S.~L. Brunton, J.~L. Proctor, and J.~N. Kutz.
\newblock Discovering governing equations from data by sparse identification of
  nonlinear dynamical systems.
\newblock {\em Proceedings of the National Academy of Sciences},
  113(15):3932--3937, 2016.

\bibitem{carlberg2013gnat}
K. Carlberg, C. Farhat, J. Cortial, and D. Amsallem.
\newblock {The GNAT method for nonlinear model reduction: effective
  implementation and application to computational fluid dynamics and turbulent
  flows}.
\newblock {\em Journal of Computational Physics}, 242:623--647, 2013.

\bibitem{casella2024statistical}
G. Casella and R.~L. Berger.
\newblock {\em {Statistical Inference}}.
\newblock CRC Press, second edition, 2024.

\bibitem{chattopadhyay2021towards}
A. Chattopadhyay, M. Mustafa, P. Hassanzadeh, E. Bach, and K. Kashinath.
\newblock Towards physically consistent data-driven weather forecasting:
  Integrating data assimilation with equivariance-preserving spatial
  transformers in a case study with {ERA}5.
\newblock {\em Geoscientific Model Development Discussions}, pages 1--23, 2021.

\bibitem{chattopadhyay2020deep}
A. Chattopadhyay, M. Mustafa, P. Hassanzadeh, and K. Kashinath.
\newblock Deep spatial transformers for autoregressive data-driven forecasting
  of geophysical turbulence.
\newblock In {\em Proceedings of the 10th International Conference on Climate
  Informatics}, pages 106--112, 2020.

\bibitem{chattopadhyay2020superparameterization}
A. Chattopadhyay, A. Subel, and P. Hassanzadeh.
\newblock Data-driven super-parameterization using deep learning:
  Experimentation with multiscale {L}orenz 96 systems and transfer learning.
\newblock {\em Journal of Advances in Modeling Earth Systems},
  12(11):e2020MS002084, 2020.

\bibitem{CLM17_L9D}
M.~D. Chekroun, H. Liu, and J.~C. McWilliams.
\newblock The emergence of fast oscillations in a reduced primitive equation
  model and its implications for closure theories.
\newblock {\em Computers \& Fluids}, 151:3--22, 2017.

\bibitem{CLM20}
M.~D. Chekroun, H. Liu, and J.~C. McWilliams.
\newblock {Variational approach to closure of nonlinear dynamical systems:
  Autonomous case}.
\newblock {\em Journal of Statistical Physics}, 179:1073--1160, 2020.

\bibitem{chekroun2021stochastic}
M.~D. Chekroun, H. Liu, and J.~C. McWilliams.
\newblock Stochastic rectification of fast oscillations on slow manifold
  closures.
\newblock {\em Proceedings of the National Academy of Sciences},
  118(48):e2113650118, 2021.

\bibitem{chekroun2017data}
M.~D. Chekroun and D. Kondrashov.
\newblock {Data-adaptive harmonic spectra and multilayer Stuart-Landau models}.
\newblock {\em Chaos: An Interdisciplinary Journal of Nonlinear Science},
  27(9):093110, 2017.

\bibitem{chekroun2011predicting}
M.~D. Chekroun, D. Kondrashov, and M. Ghil.
\newblock Predicting stochastic systems by noise sampling, and application to
  the {E}l {N}i{\~n}o-southern oscillation.
\newblock {\em Proceedings of the National Academy of Sciences},
  108(29):11766--11771, 2011.

\bibitem{chen2020learning}
N. Chen.
\newblock Learning nonlinear turbulent dynamics from partial observations via
  analytically solvable conditional statistics.
\newblock {\em Journal of Computational Physics}, 418:109635, 2020.

\bibitem{chen2021bamcafe}
N. Chen and Y. Li.
\newblock {BAMCAFE}: A {B}ayesian machine learning advanced forecast ensemble
  method for complex turbulent systems with partial observations.
\newblock {\em Chaos: An Interdisciplinary Journal of Nonlinear Science},
  31(11):113114, 2021.

\bibitem{chen2022conditional}
N. Chen, Y. Li, and H. Liu.
\newblock Conditional gaussian nonlinear system: A fast preconditioner and a
  cheap surrogate model for complex nonlinear systems.
\newblock {\em Chaos: An Interdisciplinary Journal of Nonlinear Science},
  32:053122, 2022.

\bibitem{chen2024minimum}
N. Chen and H. Liu.
\newblock Minimum reduced-order models via causal inference.
\newblock {\em arXiv preprint arXiv:2407.00271}, pages 1--31, 2024.

\bibitem{chen2022shock}
N. Chen, H. Liu, and F. Lu.
\newblock Shock trace prediction by reduced models for a viscous stochastic
  burgers equation.
\newblock {\em Chaos: An Interdisciplinary Journal of Nonlinear Science},
  32:043109, 2022.

\bibitem{chen2018conditional}
N. Chen and A. Majda.
\newblock Conditional {G}aussian systems for multiscale nonlinear stochastic
  systems: Prediction, state estimation and uncertainty quantification.
\newblock {\em Entropy}, 20(7):509, 2018.

\bibitem{chen2024physics}
N. Chen and D. Qi.
\newblock A physics-informed data-driven algorithm for ensemble forecast of
  complex turbulent systems.
\newblock {\em Applied Mathematics and Computation}, 466:128480, 2024.

\bibitem{chen2023causality}
N. Chen and Y. Zhang.
\newblock A causality-based learning approach for discovering the underlying
  dynamics of complex systems from partial observations with stochastic
  parameterization.
\newblock {\em Physica D: Nonlinear Phenomena}, 449:133743, 2023.

\bibitem{chorin2015discrete}
A.~J. Chorin and F. Lu.
\newblock Discrete approach to stochastic parametrization and dimension
  reduction in nonlinear dynamics.
\newblock {\em Proceedings of the National Academy of Sciences},
  112(32):9804--9809, 2015.

\bibitem{cortiella2021sparse}
A. Cortiella, K.-C. Park, and A. Doostan.
\newblock Sparse identification of nonlinear dynamical systems via reweighted
  l1-regularized least squares.
\newblock {\em Computer Methods in Applied Mechanics and Engineering},
  376:113620, 2021.

\bibitem{cover1999elements}
T. Cover and J. Thomas.
\newblock {\em Elements of Information Theory}.
\newblock John Wiley \& Sons, 2nd edition, 2006.

\bibitem{crommelin2004strategies}
D. Crommelin and A. Majda.
\newblock {Strategies for model reduction: Comparing different optimal bases}.
\newblock {\em J. Atmos. Sci.}, 61(17):2206--2217, 2004.

\bibitem{darbellay1999estimation}
G.~A. Darbellay and I. Vajda.
\newblock Estimation of the information by an adaptive partitioning of the
  observation space.
\newblock {\em IEEE Transactions on Information Theory}, 45:1315--1321, 1999.

\bibitem{elinger2021information}
J. Elinger.
\newblock {\em Information Theoretic Causality Measures For Parameter
  Estimation and System Identification.}
\newblock PhD thesis, Georgia Institute of Technology, Atlanta, GA, USA, 2021.

\bibitem{elinger2021causation}
J. Elinger and J. Rogers.
\newblock Causation entropy method for covariate selection in dynamic models.
\newblock In {\em 2021 American Control Conference (ACC)}, pages 2842--2847.
  IEEE, 2021.

\bibitem{fish2021entropic}
J. Fish, A. DeWitt, A.~A.~R. AlMomani, P.~J. Laurienti, and E. Bollt.
\newblock Entropic regression with neurologically motivated applications.
\newblock {\em Chaos: An Interdisciplinary Journal of Nonlinear Science},
  31(11), 2021.

\bibitem{ghil2012topics}
M. Ghil and S. Childress.
\newblock {\em Topics in geophysical fluid dynamics: atmospheric dynamics,
  dynamo theory, and climate dynamics}.
\newblock Springer Science \& Business Media, 2012.

\bibitem{hannachi2007empirical}
A. Hannachi, I.~T. Jolliffe, and D.~B. Stephenson.
\newblock {Empirical orthogonal functions and related techniques in atmospheric
  science: A review}.
\newblock {\em Int.~J.~Climatol.}, 27:1119--1152, 2007.

\bibitem{harlim2014ensemble}
J. Harlim, A. Mahdi, and A.~J. Majda.
\newblock An ensemble {K}alman filter for statistical estimation of physics
  constrained nonlinear regression models.
\newblock {\em Journal of Computational Physics}, 257:782--812, 2014.

\bibitem{hasselmann1988pips}
K. Hasselmann.
\newblock {PIPs and POPs: The reduction of complex dynamical systems using
  principal interaction and oscillation patterns}.
\newblock {\em Journal of Geophysical Research: Atmospheres},
  93(D9):11015--11021, 1988.

\bibitem{herawati2018regularized}
N. Herawati, K. Nisa, E. Setiawan, Nusyirwan, and Tiryono.
\newblock Regularized multiple regression methods to deal with severe
  multicollinearity.
\newblock {\em International Journal of Statistics and Applications},
  8:167--172, 2018.

\bibitem{hijazi2020data}
S. Hijazi, G. Stabile, A. Mola, and G. Rozza.
\newblock Data-driven pod-galerkin reduced order model for turbulent flows.
\newblock {\em Journal of Computational Physics}, 416:109513, 2020.

\bibitem{HLB96}
P. Holmes, J.~L. Lumley, and G. Berkooz.
\newblock {\em Turbulence, Coherent Structures, Dynamical Systems and
  Symmetry}.
\newblock Cambridge, 1996.

\bibitem{Holmes_al12}
P. Holmes, J.~L. Lumley, G. Berkooz, and C.~W. Rowley.
\newblock {\em {Turbulence, Coherent Structures, Dynamical Systems and
  Symmetry}}.
\newblock Cambridge University Press, Cambridge, second edition, 2012.

\bibitem{hyman1986order}
J.~M. Hyman, B. Nicolaenko, and S. Zaleski.
\newblock {Order and complexity in the Kuramoto-Sivashinsky model of weakly
  turbulent interfaces}.
\newblock {\em Physica D}, 23:265--292, 1986.

\bibitem{jardak2010comparison}
M. Jardak, I. Navon, and M. Zupanski.
\newblock {Comparison of sequential data assimilation methods for the
  Kuramoto--Sivashinsky equation}.
\newblock {\em International Journal for Numerical Methods in Fluids},
  62:374--402, 2010.

\bibitem{kaiser1977data}
J. Kaiser and W. Reed.
\newblock Data smoothing using low-pass digital filters.
\newblock {\em Review of Scientific Instruments}, 48:1447--1457, 1977.

\bibitem{kassam2005fourth}
A. Kassam and L.~N. Trefethen.
\newblock {Fourth-order time-stepping for stiff PDEs}.
\newblock {\em SIAM J.~Sci. Comp.}, 26(4):1214--1233, 2005.

\bibitem{kim2017causation}
P. Kim, J. Rogers, J. Sun, and E. Bollt.
\newblock Causation entropy identifies sparsity structure for parameter
  estimation of dynamic systems.
\newblock {\em Journal of Computational and Nonlinear Dynamics}, 12(1):011008,
  2017.

\bibitem{kleeman2011information}
R. Kleeman.
\newblock Information theory and dynamical system predictability.
\newblock {\em Entropy}, 13(3):612--649, 2011.

\bibitem{koc2022verifiability}
B. Koc, C. Mou, H. Liu, Z. Wang, G. Rozza, and T. Iliescu.
\newblock Verifiability of the data-driven variational multiscale reduced order
  model.
\newblock {\em Journal of Scientific Computing}, 93:54:1--26, 2022.

\bibitem{kondrashov2015data}
D. Kondrashov, M.~D. Chekroun, and M. Ghil.
\newblock {Data-driven non-Markovian closure models}.
\newblock {\em Physica D: Nonlinear Phenomena}, 297:33--55, 2015.

\bibitem{kozachenko1987sample}
L.~F. Kozachenko and N.~N. Leonenko.
\newblock Sample estimate of the entropy of a random vector.
\newblock {\em Problems of Information Transmission}, 23:95--102, 1987.

\bibitem{kraskov2004estimating}
A. Kraskov, H. St{\"o}gbauer, and P. Grassberger.
\newblock Estimating mutual information.
\newblock {\em Physical Review E}, 69:066138, 2004.

\bibitem{kravtsov2005multilevel}
S. Kravtsov, D. Kondrashov, and M. Ghil.
\newblock {Multilevel regression modeling of nonlinear processes: Derivation
  and applications to climatic variability}.
\newblock {\em Journal of Climate}, 18(21):4404--4424, 2005.

\bibitem{KV01}
K. Kunisch and S. Volkwein.
\newblock Galerkin proper orthogonal decomposition methods for parabolic
  problems.
\newblock {\em Numer. Math.}, 90:117--148, 2001.

\bibitem{kuramoto1976persistent}
Y. Kuramoto and T. Tsuzuki.
\newblock Persistent propagation of concentration waves in dissipative media
  far from thermal equilibrium.
\newblock {\em Prog. Theor. Phys.}, 55(2):356--369, 1976.

\bibitem{kwasniok1996reduction}
F. Kwasniok.
\newblock The reduction of complex dynamical systems using principal
  interaction patterns.
\newblock {\em Physica D: Nonlinear Phenomena}, 92(1-2):28--60, 1996.

\bibitem{kwasniok1997optimal}
F. Kwasniok.
\newblock {Optimal Galerkin approximations of partial differential equations
  using principal interaction patterns}.
\newblock {\em Physical Rev. E}, 55(5):5365, 1997.

\bibitem{laquey1975nonlinear}
R. LaQuey, S. Mahajan, P. Rutherford, and W. Tang.
\newblock Nonlinear saturation of the trapped-ion mode.
\newblock {\em Physical Review Letters}, 34:391--394, 1975.

\bibitem{larios2024nonlinear}
A. Larios and Y. Pei.
\newblock Nonlinear continuous data assimilation.
\newblock {\em Evolution Equations and Control Theory}, 13:329--348, 2024.

\bibitem{lee1998independent}
T.-W. Lee.
\newblock {\em {Independent Component Analysis: Theory and Applications}}.
\newblock Springer, 1998.

\bibitem{lin2021data}
K.~K. Lin and F. Lu.
\newblock Data-driven model reduction, wiener projections, and the
  {K}oopman-{M}ori-{Z}wanzig formalism.
\newblock {\em Journal of Computational Physics}, 424:109864, 2021.

\bibitem{lozano2022information}
A. Lozano-Dur{\'a}n and G. Arranz.
\newblock Information-theoretic formulation of dynamical systems: causality,
  modeling, and control.
\newblock {\em Physical Review Research}, 4:023195, 2022.

\bibitem{lu2017data}
F. Lu, K.~K. Lin, and A.~J. Chorin.
\newblock {Data-based stochastic model reduction for the Kuramoto--Sivashinsky
  equation}.
\newblock {\em Physica D: Nonlinear Phenomena}, 340:46--57, 2017.

\bibitem{lunasin2017evolution}
E. Lunasin and E.~S. Titi.
\newblock Finite determining parameters feedback control for distributed
  nonlinear dissipative systems -- a computational study.
\newblock {\em Evolution Equations and Control Theory}, 6:535--557, 2017.

\bibitem{majda2018model}
A.~J. Majda and N. Chen.
\newblock Model error, information barriers, state estimation and prediction in
  complex multiscale systems.
\newblock {\em Entropy}, 20(9):644, 2018.

\bibitem{majda2012physics}
A.~J. Majda and J. Harlim.
\newblock Physics constrained nonlinear regression models for time series.
\newblock {\em Nonlinearity}, 26(1):201, 2012.

\bibitem{moosavi2015efficient}
A. Moosavi, R. Stefanescu, and A. Sandu.
\newblock Efficient construction of local parametric reduced order models using
  machine learning techniques.
\newblock {\em arXiv preprint arXiv:1511.02909}, 2015.

\bibitem{mou2021data}
C. Mou, B. Koc, O. San, L.~G. Rebholz, and T. Iliescu.
\newblock Data-driven variational multiscale reduced order models.
\newblock {\em Computer Methods in Applied Mechanics and Engineering},
  373:113470, 2021.

\bibitem{mou2023combining}
C. Mou, L.~M. Smith, and N. Chen.
\newblock Combining stochastic parameterized reduced-order models with machine
  learning for data assimilation and uncertainty quantification with partial
  observations.
\newblock {\em Journal of Advances in Modeling Earth Systems},
  15(10):e2022MS003597, 2023.

\bibitem{noack2011reduced}
B.~R. Noack, M. Morzynski, and G. Tadmor.
\newblock {\em Reduced-order modelling for flow control}, volume 528.
\newblock Springer Science \& Business Media, 2011.

\bibitem{otto2019linearly}
S.~E. Otto and C.~W. Rowley.
\newblock Linearly recurrent autoencoder networks for learning dynamics.
\newblock {\em SIAM Journal on Applied Dynamical Systems}, 18:558--593, 2019.

\bibitem{pace2009sampling}
R.~K. Pace and J.~P. LeSage.
\newblock A sampling approach to estimate the log determinant used in spatial
  likelihood problems.
\newblock {\em Journal of Geographical Systems}, 11(3):209--225, 2009.

\bibitem{pawar2020data}
S. Pawar, S.~E. Ahmed, O. San, and A. Rasheed.
\newblock Data-driven recovery of hidden physics in reduced order modeling of
  fluid flows.
\newblock {\em Physics of Fluids}, 32(3):036602, 2020.

\bibitem{peherstorfer2015dynamic}
B. Peherstorfer and K. Willcox.
\newblock Dynamic data-driven reduced-order models.
\newblock {\em Computer Methods in Applied Mechanics and Engineering},
  291:21--41, 2015.

\bibitem{penland1993prediction}
C. Penland and T. Magorian.
\newblock {Prediction of Ni\~no 3 sea surface temperatures using linear inverse
  modeling}.
\newblock {\em Journal of Climate}, 6:1067--1076, 1993.

\bibitem{weigend1994time}
D.~R. Rigney, A.~L. Goldberger, W.~C. Ocasio, Y. Ichimaru, G.~B. Moody, and
  R.~G. Mark.
\newblock Multi-channel physiological data: Description and analysis ({Data Set
  B}).
\newblock In A.~S. Weigend and N.~A. Gershenfeld, editors, {\em {Time Series
  Prediction: Forecasting the Future and Understanding the Past}}, pages
  105--130. Routledge, Taylor \& Francis Group, New York, London, 1994.

\bibitem{rish2014sparse}
I. Rish and G.~Y. Grabarnik.
\newblock {\em {Sparse Modeling: Theory, Algorithms, and Applications}}.
\newblock CRC press, 2014.

\bibitem{rowley2009spectral}
C.~W. Rowley, I. Mezi{\'c}, S. Bagheri, P. Schlatter, and D.~S. Henningson.
\newblock Spectral analysis of nonlinear flows.
\newblock {\em Journal of Fluid Mechanics}, 641:115--127, 2009.

\bibitem{rudin1992nonlinear}
L.~I. Rudin, S. Osher, and E. Fatemi.
\newblock Nonlinear total variation based noise removal algorithms.
\newblock {\em Physica D}, 60:259--268, 1992.

\bibitem{san2018extreme}
O. San and R. Maulik.
\newblock Extreme learning machine for reduced order modeling of turbulent
  geophysical flows.
\newblock {\em Physical Review E}, 97(4):042322, 2018.

\bibitem{santosa1986linear}
F. Santosa and W.~W. Symes.
\newblock Linear inversion of band-limited reflection seismograms.
\newblock {\em SIAM Journal on Scientific and Statistical Computing},
  7(4):1307--1330, 1986.

\bibitem{schaeffer2018extracting}
H. Schaeffer, G. Tran, and R. Ward.
\newblock Extracting sparse high-dimensional dynamics from limited data.
\newblock {\em SIAM Journal on Applied Mathematics}, 78(6):3279--3295, 2018.

\bibitem{schmid2010dynamic}
P.~J. Schmid.
\newblock Dynamic mode decomposition of numerical and experimental data.
\newblock {\em Journal of Fluid Mechanics}, 656:5--28, 2010.

\bibitem{schneider2021learning}
T. Schneider, A.~M. Stuart, and J.-L. Wu.
\newblock Learning stochastic closures using ensemble {K}alman inversion.
\newblock {\em Transactions of Mathematics and Its Applications}, 5(1):tnab003,
  2021.

\bibitem{schreiber2000measuring}
T. Schreiber.
\newblock Measuring information transfer.
\newblock {\em Physical Review Letters}, 85:461--464, 2000.

\bibitem{seber2003linear}
G.~A.~F. Seber and A.~J. Lee.
\newblock {\em {Linear Regression Analysis}}.
\newblock John Wiley \& Sons, second edition, 2003.

\bibitem{sheard2009principles}
S.~A. Sheard and A. Mostashari.
\newblock Principles of complex systems for systems engineering.
\newblock {\em Systems Engineering}, 12(4):295--311, 2009.

\bibitem{Sir87abc}
L. Sirovich.
\newblock Turbulence and the dynamics of coherent structures. {P}arts
  {I}--{III}.
\newblock {\em Quart. Appl. Math.}, 45(3):561--590, 1987.

\bibitem{sivashinsky1977nonlinear}
G. Sivashinsky.
\newblock {Nonlinear analysis of hydrodynamic instability in laminar
  flames-{I}. Derivation of basic equations}.
\newblock {\em Acta Astronautica}, 4(11-12):1177--1206, 1977.

\bibitem{sivashinsky1980irregular}
G.~I. Sivashinsky and D.~M. Michelson.
\newblock On irregular wavy flow of a liquid film down a vertical plane.
\newblock {\em Progress of theoretical physics}, 63:2112--2114, 1980.

\bibitem{smarra2018data}
F. Smarra, A. Jain, T. De~Rubeis, D. Ambrosini, A. D’Innocenzo, and R.
  Mangharam.
\newblock Data-driven model predictive control using random forests for
  building energy optimization and climate control.
\newblock {\em Applied energy}, 226:1252--1272, 2018.

\bibitem{snyder2022reduced}
W. Snyder, C. Mou, H. Liu, O. San, R. De~Vita, and T. Iliescu.
\newblock Reduced order model closures: A brief tutorial.
\newblock In {\em Recent Advances in Mechanics and Fluid-Structure Interaction
  with Applications: The Bong Jae Chung Memorial Volume}, pages 167--193.
  Springer, 2022.

\bibitem{srinivasan2024turbulence}
K. Srinivasan, M.~D. Chekroun, and J.~C. McWilliams.
\newblock {Turbulence closure with small, local neural networks: Forced
  two-dimensional and $\beta$-plane flows}.
\newblock {\em Journal of Advances in Modeling Earth Systems},
  16:e2023MS003795, 2024.

\bibitem{stinis2004stochastic}
P. Stinis.
\newblock {Stochastic optimal prediction for the Kuramoto--Sivashinsky
  equation}.
\newblock {\em Multiscale Model. Simul.}, 2(4):580--612, 2004.

\bibitem{strogatz2018nonlinear}
S.~H. Strogatz.
\newblock {\em Nonlinear dynamics and chaos with student solutions manual: With
  applications to physics, biology, chemistry, and engineering}.
\newblock CRC press, 2018.

\bibitem{sun2014causation}
J. Sun and E.~M. Bollt.
\newblock Causation entropy identifies indirect influences, dominance of
  neighbors and anticipatory couplings.
\newblock {\em Physica D}, 267:49--57, 2014.

\bibitem{sun2015causal}
J. Sun, D. Taylor, and E.~M. Bollt.
\newblock Causal network inference by optimal causation entropy.
\newblock {\em SIAM Journal on Applied Dynamical Systems}, 14:73--106, 2015.

\bibitem{taira2017modal}
K. Taira, S.~L. Brunton, S.~T.~M. Dawson, C.~W. Rowley, T. Colonius, B.~J.
  McKeon, O.~T. Schmidt, S. Gordeyev, V. Theofilis, and L.~S. Ukeiley.
\newblock {Modal analysis of fluid flows: An overview}.
\newblock {\em AIAA Journal}, 55:4013--4041, 2017.

\bibitem{taira2020modal}
K. Taira, M.~S. Hemati, S.~L. Brunton, Y. Sun, K. Duraisamy, S. Bagheri, S.~T.
  Dawson, and C.-A. Yeh.
\newblock Modal analysis of fluid flows: Applications and outlook.
\newblock {\em AIAA journal}, 58(3):998--1022, 2020.

\bibitem{temam1997infinite}
R. Temam.
\newblock {\em {Infinite-Dimensional Dynamical Systems in Mechanics and
  Physics}}, volume~68 of {\em Applied Mathematical Sciences}.
\newblock Springer, New York, 2nd edition, 1997.

\bibitem{tibshirani1996regression}
R. Tibshirani.
\newblock Regression shrinkage and selection via the lasso.
\newblock {\em Journal of the Royal Statistical Society: Series B
  (Methodological)}, 58(1):267--288, 1996.

\bibitem{tippett2004measuring}
M.~K. Tippett, R. Kleeman, and Y. Tang.
\newblock Measuring the potential utility of seasonal climate predictions.
\newblock {\em Geophysical research letters}, 31(22), 2004.

\bibitem{tu2013dynamic}
J.~H. Tu, C.~W. Rowley, D.~M. Luchtenburg, S.~L. Brunton, and J.~N. Kutz.
\newblock {On dynamic mode decomposition: Theory and applications}.
\newblock {\em Journal of Computational Dynamics}, 1:391--421, 2014.

\bibitem{vallis2017atmospheric}
G.~K. Vallis.
\newblock {\em Atmospheric and oceanic fluid dynamics}.
\newblock Cambridge University Press, 2017.

\bibitem{vautard1992singular}
R. Vautard, P. Yiou, and M. Ghil.
\newblock {Singular-spectrum analysis: A toolkit for short, noisy chaotic
  signals}.
\newblock {\em Physica D}, 58:95--126, 1992.

\bibitem{wan2017reduced}
Z.~Y. Wan and T.~P. Sapsis.
\newblock Reduced-space {G}aussian process regression for data-driven
  probabilistic forecast of chaotic dynamical systems.
\newblock {\em Physica D: Nonlinear Phenomena}, 345:40--55, 2017.

\bibitem{wilcox1988multiscale}
D.~C. Wilcox.
\newblock Multiscale model for turbulent flows.
\newblock {\em AIAA journal}, 26(11):1311--1320, 1988.

\bibitem{williams2015data}
M. Williams, I. Kevrekidis, and C. Rowley.
\newblock {A data--driven approximation of the Koopman operator: Extending
  dynamic mode decomposition}.
\newblock {\em Journal of Nonlinear Science}, 25(6):1307--1346, 2015.

\bibitem{wolf1985determining}
A. Wolf, J.~B. Swift, H.~L. Swinney, and J.~A. Vastano.
\newblock Determining lyapunov exponents from a time series.
\newblock {\em Physica D}, 16:285--317, 1985.

\bibitem{wyner1978definition}
A.~D. Wyner.
\newblock A definition of conditional mutual information for arbitrary
  ensembles.
\newblock {\em Information and Control}, 38:51--59, 1978.

\bibitem{xie2018data}
X. Xie, M. Mohebujjaman, L.~G. Rebholz, and T. Iliescu.
\newblock Data-driven filtered reduced order modeling of fluid flows.
\newblock {\em SIAM Journal on Scientific Computing}, 40(3):B834--B857, 2018.

\end{thebibliography}
\end{document}